

\input amstex
\expandafter\ifx\csname mathdefs.tex\endcsname\relax
  \expandafter\gdef\csname mathdefs.tex\endcsname{}
\else \message{Hey!  Apparently you were trying to
  \string twice.   This does not make sense.} 
\errmessage{Please edit your file (probably \jobname.tex) and remove
any duplicate ``\string\input'' lines} \fi




\catcode`\X=12\catcode`\@=11

\def\n@wcount{\alloc@0\count\countdef\insc@unt}
\def\n@wwrite{\alloc@7\write\chardef\sixt@@n}
\def\n@wread{\alloc@6\read\chardef\sixt@@n}
\def\r@s@t{\relax}\def\v@idline{\par}\def\@mputate#1/{#1}
\def\l@c@l#1X{\firstpart.#1}\def\gl@b@l#1X{#1}\def\t@d@l#1X{{}}

\def\crossrefs#1{\ifx\all#1\let\tr@ce=\all\else\def\tr@ce{#1,}\fi
   \n@wwrite\cit@tionsout\openout\cit@tionsout=\jobname.cit 
   \write\cit@tionsout{\tr@ce}\expandafter\setfl@gs\tr@ce,}
\def\setfl@gs#1,{\def\@{#1}\ifx\@\empty\let\next=\relax
   \else\let\next=\setfl@gs\expandafter\xdef
   \csname#1tr@cetrue\endcsname{}\fi\next}
\def\m@ketag#1#2{\expandafter\n@wcount\csname#2tagno\endcsname
     \csname#2tagno\endcsname=0\let\tail=\all\xdef\all{\tail#2,}
   \ifx#1\l@c@l\let\tail=\r@s@t\xdef\r@s@t{\csname#2tagno\endcsname=0\tail}\fi
   \expandafter\gdef\csname#2cite\endcsname##1{\expandafter
     \ifx\csname#2tag##1\endcsname\relax?\else\csname#2tag##1\endcsname\fi
     \expandafter\ifx\csname#2tr@cetrue\endcsname\relax\else
     \write\cit@tionsout{#2tag ##1 cited on page \folio.}\fi}
   \expandafter\gdef\csname#2page\endcsname##1{\expandafter
     \ifx\csname#2page##1\endcsname\relax?\else\csname#2page##1\endcsname\fi
     \expandafter\ifx\csname#2tr@cetrue\endcsname\relax\else
     \write\cit@tionsout{#2tag ##1 cited on page \folio.}\fi}
   \expandafter\gdef\csname#2tag\endcsname##1{\expandafter
      \ifx\csname#2check##1\endcsname\relax
      \expandafter\xdef\csname#2check##1\endcsname{}%
      \else\immediate\write16{Warning: #2tag ##1 used more than once.}\fi
      \multit@g{#1}{#2}##1/X%
      \write\t@gsout{#2tag ##1 assigned number \csname#2tag##1\endcsname\space
      on page \number\count0.}%
   \csname#2tag##1\endcsname}}
\def\multit@g#1#2#3/#4X{\def\t@mp{#4}\ifx\t@mp\empty%
      \global\advance\csname#2tagno\endcsname by 1 
      \expandafter\xdef\csname#2tag#3\endcsname
      {#1\number\csname#2tagno\endcsnameX}%
   \else\expandafter\ifx\csname#2last#3\endcsname\relax
      \expandafter\n@wcount\csname#2last#3\endcsname
      \global\advance\csname#2tagno\endcsname by 1 
      \expandafter\xdef\csname#2tag#3\endcsname
      {#1\number\csname#2tagno\endcsnameX}
      \write\t@gsout{#2tag #3 assigned number \csname#2tag#3\endcsname\space
      on page \number\count0.}\fi
   \global\advance\csname#2last#3\endcsname by 1
   \def\t@mp{\expandafter\xdef\csname#2tag#3/}%
   \expandafter\t@mp\@mputate#4\endcsname
   {\csname#2tag#3\endcsname\lastpart{\csname#2last#3\endcsname}}\fi}
\def\t@gs#1{\def\all{}\m@ketag#1e\m@ketag#1s\m@ketag\t@d@l p
   \m@ketag\gl@b@l r \n@wread\t@gsin
   \openin\t@gsin=\jobname.tgs \re@der \closein\t@gsin
   \n@wwrite\t@gsout\openout\t@gsout=\jobname.tgs }
\outer\def\localtags{\t@gs\l@c@l}
\outer\def\globaltags{\t@gs\gl@b@l}
\outer\def\newlocaltag#1{\m@ketag\l@c@l{#1}}
\outer\def\newglobaltag#1{\m@ketag\gl@b@l{#1}}

\newif\ifpr@ 
\def\m@kecs #1tag #2 assigned number #3 on page #4.%
   {\expandafter\gdef\csname#1tag#2\endcsname{#3}
   \expandafter\gdef\csname#1page#2\endcsname{#4}
   \ifpr@\expandafter\xdef\csname#1check#2\endcsname{}\fi}
\def\re@der{\ifeof\t@gsin\let\next=\relax\else
   \read\t@gsin to\t@gline\ifx\t@gline\v@idline\else
   \expandafter\m@kecs \t@gline\fi\let \next=\re@der\fi\next}
\def\pretags#1{\pr@true\pret@gs#1,,}
\def\pret@gs#1,{\def\@{#1}\ifx\@\empty\let\n@xtfile=\relax
   \else\let\n@xtfile=\pret@gs \openin\t@gsin=#1.tgs \message{#1} \re@der 
   \closein\t@gsin\fi \n@xtfile}

\newcount\sectno\sectno=0\newcount\subsectno\subsectno=0
\newif\ifultr@local \def\ultralocal{\ultr@localtrue}
\def\firstpart{\number\sectno}
\def\lastpart#1{\ifcase#1 \or a\or b\or c\or d\or e\or f\or g\or h\or 
   i\or k\or l\or m\or n\or o\or p\or q\or r\or s\or t\or u\or v\or w\or 
   x\or y\or z \fi}

\def\resetall{\global\advance\sectno by 1\subsectno=0
   \gdef\firstpart{\number\sectno}\r@s@t}
\def\resetsub{\global\advance\subsectno by 1
   \gdef\firstpart{\number\sectno.\number\subsectno}\r@s@t}
\def\newsection#1\par{\resetall\vskip0pt plus.3\vsize\penalty-250
   \vskip0pt plus-.3\vsize\bigskip\bigskip
   \message{#1}\leftline{\bf#1}\nobreak\bigskip}
\def\subsection#1\par{\ifultr@local\resetsub\fi
   \vskip0pt plus.2\vsize\penalty-250\vskip0pt plus-.2\vsize
   \bigskip\smallskip\message{#1}\leftline{\bf#1}\nobreak\medskip}

\def\t@gsoff#1,{\def\@{#1}\ifx\@\empty\let\next=\relax\else\let\next=\t@gsoff
   \def\@@{p}\ifx\@\@@\else
   \expandafter\gdef\csname#1cite\endcsname##1{\zeigen{##1}}
   \expandafter\gdef\csname#1page\endcsname##1{?}
   \expandafter\gdef\csname#1tag\endcsname##1{\zeigen{##1}}\fi\fi\next}
\def\verbatimtags{\ifx\all\relax\else\expandafter\t@gsoff\all,\fi}
\def\zeigen#1{\hbox{$\langle$}#1\hbox{$\rangle$}}

\def\(#1){\edef\dot@g{\ifmmode\ifinner(\hbox{\noexpand\etag{#1}})
   \else\noexpand\eqno(\hbox{\noexpand\etag{#1}})\fi
   \else(\noexpand\ecite{#1})\fi}\dot@g}

\newif\ifbr@ck
\def\eat#1{}
\def\[#1]{\br@cktrue[\br@cket#1'X]}
\def\br@cket#1'#2X{\def\temp{#2}\ifx\temp\empty\let\next\eat
   \else\let\next\br@cket\fi
   \ifbr@ck\br@ckfalse\br@ck@t#1,X\else\br@cktrue#1\fi\next#2X}
\def\br@ck@t#1,#2X{\def\temp{#2}\ifx\temp\empty\let\neext\eat
   \else\let\neext\br@ck@t\def\temp{,}\fi
   \def\teemp{#1}\ifx\teemp\empty\else\rcite{#1}\fi\temp\neext#2X}
\def\resetbr@cket{\gdef\[##1]{[\rtag{##1}]}}
\def\references{\resetbr@cket\newsection References\par}

\newtoks\symb@ls\newtoks\s@mb@ls\newtoks\p@gelist\n@wcount\ftn@mber
    \ftn@mber=1\newif\ifftn@mbers\ftn@mbersfalse\newif\ifbyp@ge\byp@gefalse
\def\defm@rk{\ifftn@mbers\n@mberm@rk\else\symb@lm@rk\fi}
\def\n@mberm@rk{\xdef\m@rk{{\the\ftn@mber}}%
    \global\advance\ftn@mber by 1 }
\def\rot@te#1{\let\temp=#1\global#1=\expandafter\r@t@te\the\temp,X}
\def\r@t@te#1,#2X{{#2#1}\xdef\m@rk{{#1}}}
\def\b@@st#1{{$^{#1}$}}\def\str@p#1{#1}
\def\symb@lm@rk{\ifbyp@ge\rot@te\p@gelist\ifnum\expandafter\str@p\m@rk=1 
    \s@mb@ls=\symb@ls\fi\write\f@nsout{\number\count0}\fi \rot@te\s@mb@ls}
\def\byp@ge{\byp@getrue\n@wwrite\f@nsin\openin\f@nsin=\jobname.fns 
    \n@wcount\currentp@ge\currentp@ge=0\p@gelist={0}
    \re@dfns\closein\f@nsin\rot@te\p@gelist
    \n@wread\f@nsout\openout\f@nsout=\jobname.fns }
\def\m@kelist#1X#2{{#1,#2}}
\def\re@dfns{\ifeof\f@nsin\let\next=\relax\else\read\f@nsin to \f@nline
    \ifx\f@nline\v@idline\else\let\t@mplist=\p@gelist
    \ifnum\currentp@ge=\f@nline
    \global\p@gelist=\expandafter\m@kelist\the\t@mplistX0
    \else\currentp@ge=\f@nline
    \global\p@gelist=\expandafter\m@kelist\the\t@mplistX1\fi\fi
    \let\next=\re@dfns\fi\next}
\def\symbols#1{\symb@ls={#1}\s@mb@ls=\symb@ls} 
\def\bigsymbol{\textstyle}
\symbols{\bigsymbol\ast,\dagger,\ddagger,\sharp,\flat,\natural,\star}
\def\ftnumbers{\ftn@mberstrue} \def\ftsymbols{\ftn@mbersfalse}
\def\paginal{\byp@ge} \def\resetftnumbers{\ftn@mber=1}
\def\ftnote#1{\defm@rk\expandafter\expandafter\expandafter\footnote
    \expandafter\b@@st\m@rk{#1}}

\long\def\jump#1\endjump{}
\def\ssum{\mathop{\lower .1em\hbox{$\textstyle\Sigma$}}\nolimits}

\def\qed{\nobreak\kern 1em \vrule height .5em width .5em depth 0em}
\def\newneq{\hbox{\rlap{\hbox to 1\wd9{\hss$=$\hss}}\raise .1em 
   \hbox to 1\wd9{\hss$\scriptscriptstyle/$\hss}}}
\def\subsetne{\setbox9 = \hbox{$\subset$}\mathrel{\hbox{\rlap
   {\lower .4em \newneq}\raise .13em \hbox{$\subset$}}}}
\def\supsetne{\setbox9 = \hbox{$\subset$}\mathrel{\hbox{\rlap
   {\lower .4em \newneq}\raise .13em \hbox{$\supset$}}}}

\def\vbar{\mathchoice{\vrule height6.3ptdepth-.5ptwidth.8pt\kern-.8pt}
   {\vrule height6.3ptdepth-.5ptwidth.8pt\kern-.8pt}
   {\vrule height4.1ptdepth-.35ptwidth.6pt\kern-.6pt}
   {\vrule height3.1ptdepth-.25ptwidth.5pt\kern-.5pt}}
\def\f@dge{\mathchoice{}{}{\mkern.5mu}{\mkern.8mu}}
\def\b@c#1#2{{\rm \mkern#2mu\vbar\mkern-#2mu#1}}
\def\b@b#1{{\rm I\mkern-3.5mu #1}}
\def\b@a#1#2{{\rm #1\mkern-#2mu\f@dge #1}}
\def\bb#1{{\count4=`#1 \advance\count4by-64 \ifcase\count4\or\b@a A{11.5}\or
   \b@b B\or\b@c C{5}\or\b@b D\or\b@b E\or\b@b F \or\b@c G{5}\or\b@b H\or
   \b@b I\or\b@c J{3}\or\b@b K\or\b@b L \or\b@b M\or\b@b N\or\b@c O{5} \or
   \b@b P\or\b@c Q{5}\or\b@b R\or\b@a S{8}\or\b@a T{10.5}\or\b@c U{5}\or
   \b@a V{12}\or\b@a W{16.5}\or\b@a X{11}\or\b@a Y{11.7}\or\b@a Z{7.5}\fi}}

\catcode`\X=11 \catcode`\@=12

\expandafter\ifx\csname citeadd.tex\endcsname\relax
\expandafter\gdef\csname citeadd.tex\endcsname{}
\else \message{Hey!  Apparently you were trying to
\string twice.   This does not make sense.} 
\errmessage{Please edit your file (probably \jobname.tex) and remove
any duplicate ``\string\input'' lines} \fi

\sectno=-1   
\localtags
\NoBlackBoxes
\ifx\shlhetal\undefinedcontrolsequence\let\shlhetal\relax\fi
\define\mr{\medskip\roster}
\define\sn{\smallskip\noindent}
\define\mn{\medskip\noindent}
\define\bn{\bigskip\noindent}
\define\ub{\underbar}
\define\wilog{\text{without loss of generality}}
\define\ermn{\endroster\medskip\noindent}
\define\dbca{\dsize\bigcap}
\define\dbcu{\dsize\bigcup}
\define \nl{\newline}
\documentstyle {amsppt}
\topmatter
\title{Not collapsing cardinal $\le \kappa$ in $(< \kappa)$-support
iteration \\
 Part II \\
Sh667} \endtitle
\rightheadtext{Not collapsing cardinal}
\author {Saharon Shelah \thanks {\null\newline I would like to thank 
Alice Leonhardt for the beautiful typing. \null\newline
 Notes - Spring '96 \null\newline
 First Typed - 97/June/30 \null\newline
Latest Revision - 98/Aug/26} \endthanks} \endauthor 
\affil{Institute of Mathematics\\
 The Hebrew University\\
 Jerusalem, Israel
 \medskip
 Rutgers University\\
 Mathematics Department\\
 New Brunswick, NJ  USA} \endaffil

\abstract {We deal with $(< \kappa)$-supported iterated forcing notions which
are $(\hat{\Cal E}_0,\hat{\Cal E}_1)$-complete, have in mind problems on
Whitehead groups, uniformizations and the general problem.  We deal mainly
with the successor of a singular case.  This continues \cite{Sh:587}.
We also deal with complimentary
combinatorial results.} \endabstract
\endtopmatter
\document  

\expandafter\ifx\csname alice2jlem.tex\endcsname\relax
  \expandafter\gdef\csname alice2jlem.tex\endcsname{}
\else \message{Hey!  Apparently you were trying to
\string  twice.   This does not make sense.}
\errmessage{Please edit your file (probably \jobname.tex) and remove
any duplicate ``\string\input'' lines} \fi

\expandafter\ifx\csname bib4plain.tex\endcsname\relax
  \expandafter\gdef\csname bib4plain.tex\endcsname{}
\else \message{Hey!  Apparently you were trying to \string twice.   This does not make sense.}
\errmessage{Please edit your file (probably \jobname.tex) and remove
any duplicate ``\string\input'' lines} \fi

\def\renewcommand{\newcommand}	       
\edef\cite{\the\catcode`@}%
\catcode`@ = 11
\let\@oldatcatcode = \cite
\chardef\@letter = 11
\chardef\@other = 12
%
%
%
%
\def\@innerdef#1#2{\edef#1{\expandafter\noexpand\csname #2\endcsname}}%
%
%
\@innerdef\@innernewcount{newcount}%
\@innerdef\@innernewdimen{newdimen}%
\@innerdef\@innernewif{newif}%
\@innerdef\@innernewwrite{newwrite}%
%
%
%
\def\@gobble#1{}%
%
%
%
\ifx\inputlineno\@undefined
   \let\@linenumber = \empty 
\else
   \def\@linenumber{\the\inputlineno:\space}%
\fi
%
%
%
\def\@futurenonspacelet#1{\def\cs{#1}%
   \afterassignment\@stepone\let\@nexttoken=
}%
\begingroup 
\def\\{\global\let\@stoken= }%
\\ 
\endgroup
\def\@stepone{\expandafter\futurelet\cs\@steptwo}%
\def\@steptwo{\expandafter\ifx\cs\@stoken\let\@@next=\@stepthree
   \else\let\@@next=\@nexttoken\fi \@@next}%
\def\@stepthree{\afterassignment\@stepone\let\@@next= }%
%
%
%
\def\@getoptionalarg#1{%
   \let\@optionaltemp = #1%
   \let\@optionalnext = \relax
   \@futurenonspacelet\@optionalnext\@bracketcheck
}%
%
%
\def\@bracketcheck{%
   \ifx [\@optionalnext
      \expandafter\@@getoptionalarg
   \else
      \let\@optionalarg = \empty
      \expandafter\@optionaltemp
   \fi
}%
\def\@@getoptionalarg[#1]{%
   \def\@optionalarg{#1}%
   \@optionaltemp
}%
%
%
%
\def\@nnil{\@nil}%
\def\@fornoop#1\@@#2#3{}%
\def\@for#1:=#2\do#3{%
   \edef\@fortmp{#2}%
   \ifx\@fortmp\empty \else
      \expandafter\@forloop#2,\@nil,\@nil\@@#1{#3}%
   \fi
}%
\def\@forloop#1,#2,#3\@@#4#5{\def#4{#1}\ifx #4\@nnil \else
       #5\def#4{#2}\ifx #4\@nnil \else#5\@iforloop #3\@@#4{#5}\fi\fi
}%
\def\@iforloop#1,#2\@@#3#4{\def#3{#1}\ifx #3\@nnil
       \let\@nextwhile=\@fornoop \else
      #4\relax\let\@nextwhile=\@iforloop\fi\@nextwhile#2\@@#3{#4}%
}%
%
%
%
\@innernewif\if@fileexists
\def\@testfileexistence{\@getoptionalarg\@finishtestfileexistence}%
\def\@finishtestfileexistence#1{%
   \begingroup
      \def\extension{#1}%
      \immediate\openin0 =
         \ifx\@optionalarg\empty\jobname\else\@optionalarg\fi
         \ifx\extension\empty \else .#1\fi
         \space
      \ifeof 0
         \global\@fileexistsfalse
      \else
         \global\@fileexiststrue
      \fi
      \immediate\closein0
   \endgroup
}%
%
%
%
%
\def\bibliographystyle#1{%
   \@readauxfile
   \@writeaux{\string\bibstyle{#1}}%
}%
\let\bibstyle = \@gobble
%
%
\let\bblfilebasename = \jobname
\def\bibliography#1{%
   \@readauxfile
   \@writeaux{\string\bibdata{#1}}%
   \@testfileexistence[\bblfilebasename]{bbl}%
   \if@fileexists
      \nobreak
      \@readbblfile
   \fi
}%
\let\bibdata = \@gobble
%
%
\def\nocite#1{%
   \@readauxfile
   \@writeaux{\string\citation{#1}}%
}%
\@innernewif\if@notfirstcitation
%
%
\def\cite{\@getoptionalarg\@cite}%
%
%
\def\@cite#1{%
   \let\@citenotetext = \@optionalarg
   \printcitestart
   \nocite{#1}%
   \@notfirstcitationfalse
   \@for \@citation :=#1\do
   {%
      \expandafter\@onecitation\@citation\@@
   }%
   \ifx\empty\@citenotetext\else
      \printcitenote{\@citenotetext}%
   \fi
   \printcitefinish
}%
\def\@onecitation#1\@@{%
   \if@notfirstcitation
      \printbetweencitations
   \fi
   \expandafter \ifx \csname\@citelabel{#1}\endcsname \relax
      \if@citewarning
         \message{\@linenumber Undefined citation `#1'.}%
      \fi
      \expandafter\gdef\csname\@citelabel{#1}\endcsname{%
\strut
\vadjust{\vskip-\dp\strutbox
\vbox to 0pt{\vss\parindent0cm \leftskip=\hsize 
\advance\leftskip3mm
\advance\hsize 4cm\strut\openup-4pt 
\rightskip 0cm plus 1cm minus 0.5cm ?  #1 ?\strut}}
         {\tt
            \escapechar = -1
            \nobreak\hskip0pt
            \expandafter\string\csname#1\endcsname
            \nobreak\hskip0pt
         }%
      }%
   \fi
   \csname\@citelabel{#1}\endcsname
   \@notfirstcitationtrue
}%
%
%
\def\@citelabel#1{b@#1}%
%
%
\def\@citedef#1#2{\expandafter\gdef\csname\@citelabel{#1}\endcsname{#2}}%
%
%
%
\def\@readbblfile{%
   \ifx\@itemnum\@undefined
      \@innernewcount\@itemnum
   \fi
   \begingroup
      \def\begin##1##2{%
         \setbox0 = \hbox{\biblabelcontents{##2}}%
         \biblabelwidth = \wd0
      }%
      \def\end##1{}
      %
      %
      \@itemnum = 0
      \def\bibitem{\@getoptionalarg\@bibitem}%
      \def\@bibitem{%
         \ifx\@optionalarg\empty
            \expandafter\@numberedbibitem
         \else
            \expandafter\@alphabibitem
         \fi
      }%
      \def\@alphabibitem##1{%
         \expandafter \xdef\csname\@citelabel{##1}\endcsname {\@optionalarg}%
         \ifx\biblabelprecontents\@undefined
            \let\biblabelprecontents = \relax
         \fi
         \ifx\biblabelpostcontents\@undefined
            \let\biblabelpostcontents = \hss
         \fi
         \@finishbibitem{##1}%
      }%
      \def\@numberedbibitem##1{%
         \advance\@itemnum by 1
         \expandafter \xdef\csname\@citelabel{##1}\endcsname{\number\@itemnum}%
         \ifx\biblabelprecontents\@undefined
            \let\biblabelprecontents = \hss
         \fi
         \ifx\biblabelpostcontents\@undefined
            \let\biblabelpostcontents = \relax
         \fi
         \@finishbibitem{##1}%
      }%
      \def\@finishbibitem##1{%
         \biblabelprint{\csname\@citelabel{##1}\endcsname}%
         \@writeaux{\string\@citedef{##1}{\csname\@citelabel{##1}\endcsname}}%
         \ignorespaces
      }%
      %
      %
      \let\em = \bblem
      \let\newblock = \bblnewblock
      \let\sc = \bblsc
      \frenchspacing
      \clubpenalty = 4000 \widowpenalty = 4000
      \tolerance = 10000 \hfuzz = .5pt
      \everypar = {\hangindent = \biblabelwidth
                      \advance\hangindent by \biblabelextraspace}%
      \bblrm
      \parskip = 1.5ex plus .5ex minus .5ex
      \biblabelextraspace = .5em
      \bblhook
      \input \bblfilebasename.bbl
   \endgroup
}%
%
%
\@innernewdimen\biblabelwidth
\@innernewdimen\biblabelextraspace
%
%
%
\def\biblabelprint#1{%
   \noindent
   \hbox to \biblabelwidth{%
      \biblabelprecontents
      \biblabelcontents{#1}%
      \biblabelpostcontents
   }%
   \kern\biblabelextraspace
}%
%
%
%
\def\biblabelcontents#1{{\bblrm [#1]}}%
%
%
\def\bblrm{\rm}%
%
%
\def\bblem{\it}%
%
%
\def\bblsc{\ifx\@scfont\@undefined
              \font\@scfont = cmcsc10
           \fi
           \@scfont
}%
%
%
\def\bblnewblock{\hskip .11em plus .33em minus .07em }%
%
%
\let\bblhook = \empty
%
%
%
\def\printcitestart{[}
\def\printcitefinish{]}
\def\printbetweencitations{, }
\def\printcitenote#1{, #1}
%
%
%
\let\citation = \@gobble
%
%
%
\@innernewcount\@numparams
%
%
\def\newcommand#1{%
   \def\@commandname{#1}%
   \@getoptionalarg\@continuenewcommand
}%
%
%
\def\@continuenewcommand{%
   \@numparams = \ifx\@optionalarg\empty 0\else\@optionalarg \fi \relax
   \@newcommand
}%
%
%
\def\@newcommand#1{%
   \def\@startdef{\expandafter\edef\@commandname}%
   \ifnum\@numparams=0
      \let\@paramdef = \empty
   \else
      \ifnum\@numparams>9
         \errmessage{\the\@numparams\space is too many parameters}%
      \else
         \ifnum\@numparams<0
            \errmessage{\the\@numparams\space is too few parameters}%
         \else
            \edef\@paramdef{%
               \ifcase\@numparams
                  \empty  No arguments.
               \or ####1%
               \or ####1####2%
               \or ####1####2####3%
               \or ####1####2####3####4%
               \or ####1####2####3####4####5%
               \or ####1####2####3####4####5####6%
               \or ####1####2####3####4####5####6####7%
               \or ####1####2####3####4####5####6####7####8%
               \or ####1####2####3####4####5####6####7####8####9%
               \fi
            }%
         \fi
      \fi
   \fi
   \expandafter\@startdef\@paramdef{#1}%
}%
%
%
%
%
\def\@readauxfile{%
   \if@auxfiledone \else 
      \global\@auxfiledonetrue
      \@testfileexistence{aux}%
      \if@fileexists
         \begingroup
            \endlinechar = -1
            \catcode`@ = 11
            \input \jobname.aux
         \endgroup
      \else
         \message{\@undefinedmessage}%
         \global\@citewarningfalse
      \fi
      \immediate\openout\@auxfile = \jobname.aux
   \fi
}%
%
%
\newif\if@auxfiledone
\ifx\noauxfile\@undefined \else \@auxfiledonetrue\fi
%
%
%
%
\@innernewwrite\@auxfile
\def\@writeaux#1{\ifx\noauxfile\@undefined \write\@auxfile{#1}\fi}%
%
%
%
\ifx\@undefinedmessage\@undefined
   \def\@undefinedmessage{No .aux file; I won't give you warnings about
                          undefined citations.}%
\fi
%
%
\@innernewif\if@citewarning
\ifx\noauxfile\@undefined \@citewarningtrue\fi
%
%
%
\catcode`@ = \@oldatcatcode


\def\widestnumber#1#2{}

\def\rm{\fam0 \tenrm}

\def\fakesubhead#1\endsubhead{\bigskip\noindent{\bf#1}\par}


%
%
%

%

\font\textrsfs=rsfs10
\font\scriptrsfs=rsfs7
\font\scriptscriptrsfs=rsfs5

\newfam\rsfsfam
\textfont\rsfsfam=\textrsfs
\scriptfont\rsfsfam=\scriptrsfs
\scriptscriptfont\rsfsfam=\scriptscriptrsfs

\edef\oldcatcodeofat{\the\catcode`\@}
\catcode`\@11

\def\Cal@@#1{\noaccents@ \fam \rsfsfam #1}

\catcode`\@\oldcatcodeofat

\newpage

\head {Annotated Content} \endhead  \resetall 
\bn
\S1 $\quad$ GCH implies for successor of singular no stationary $S$ has
uniformization
\roster
\item "{{}}"  [For $\lambda$ strong limit singular, 
for stationary $S \subseteq S^{\lambda^+}_{\text{cf}(\lambda)}$ we 
prove strong negation of uniformization for some
$S$-ladder system and even weak versions of diamond.  E.g. if $\lambda$ is
singular strong limit $2^\lambda = \lambda^+$, then there are 
$\gamma^\delta_i < \delta$
increasing in $i < \text{ cf}(\lambda)$ with limit $\delta$ such that for
every $f:\lambda^+ \rightarrow \alpha^* < \lambda$ for stationarily many
$\delta \in S$, for every $i,f(\gamma^\delta_{2i}) = f(\gamma^\delta
_{2i+1})$.]
\endroster
\bn
\S2 $\quad$ Forcing for successor of singulars
\roster
\item "{{}}"  [Let $\lambda$ be strong limit singular $\kappa = \lambda^+
= 2^\lambda,S \subseteq S^\kappa_{\text{cf}(\lambda)}$ stationarily not
reflecting.  We present the consistency of a forcing axiom saying e.g.: if
$h_\delta$ is a function from $A_\delta$ to $\theta,A_\delta \subseteq
\delta = \sup(A_\delta)$, otp$(A_\delta) = \text{ cf}(\lambda),\theta < 
\lambda$ then for some $h:\kappa \rightarrow \theta$ for every $\delta \in
S$ we have $h_\delta \subseteq^* h$.]
\endroster
\bn
\S3 $\quad$ $\kappa^+$-c.c. and $\kappa^+$-pic
\roster
\item "{{}}"  [In the forcing axioms we would like to allow forcing notions
of cardinality $> \kappa$; for this we use a suitable chain condition
(allowed here and in \nl
\cite{Sh:587}.]
\endroster
\bn
\S4 $\quad$ Existence of non-free Whitehead groups (and Ext$(G,\Bbb Z) = 0$)
abelian groups \nl

$\quad$ in successor of singulars
\roster
\item "{{}}"  [We use the information on the existence of weak version of the
diamond for $S \subseteq S^{\lambda^+}_{\text{cf}(\lambda)},\lambda$ strong
limit singular with $2^\lambda = \lambda^+$, to prove that there are some 
abelian groups with special properties (from reasonable assumptions).]
\endroster
\newpage

\head {\S1 GCH implies for successor of singular no stationary $S$ has
unformization} \endhead  \resetall
\bn
We show that the improvement in \cite{Sh:587} over \cite{Sh:186} for
inaccessible (every ladder on $S$ rather than some) cannot be done for
successor of singulars.  This is continued in \S4. \nl
\ub{\stag{1.1} Fact}:  Assume
\mr
\item "{$(a)$}"  $\lambda$ is strong limit singular with $2^\lambda =
\lambda^+$, let cf$(\lambda) = \sigma$
\sn
\item "{$(b)$}"  $S \subseteq \{\delta < \lambda^+:\text{cf}(\delta) =
\sigma\}$ is stationary.
\ermn
\ub{Then} we can find 
$\langle < \gamma^\delta_i:i < \sigma >:\delta \in S \rangle$
such that
\mr
\item "{$(\alpha)$}"  $\gamma^\delta_i$ is increasing (in $i$) with 
limit $\delta$
\sn
\item "{$(\beta)$}"  if $\mu < \lambda$ and $f:\lambda^+ \rightarrow \mu$
\ub{then} the following set is stationary: \nl
$\{\delta \in S:f(\gamma^\delta_{2_i}) = f(\gamma^\delta_{2i+1})$ for every
$i < \sigma\}$. \nl
Moreover
\sn
\item "{$(\beta)^+$}"  if $f_i:\lambda^+ \rightarrow \mu_i,\mu_i < \lambda$
for $i < \sigma$ then the following set is stationary: \nl
$\{\delta \in S:f_i(\gamma^\delta_{2i}) = f_i(\gamma^\delta_{2i+1})$ for
every $i < \sigma\}$.
\endroster
\bigskip

\demo{Proof}  This will prove \scite{1.2}, too.  We first concentrate on
$(\alpha) + (\beta)$ only. \nl
Let $\lambda = \dsize \sum_{i < \sigma} \lambda_i,\lambda_i$ increasingly
continuous, $\lambda_{i+1} > 2^{\lambda_i},\lambda_0 > 2^\sigma$.  For
$\alpha < \lambda^+$, let $\alpha = \dbcu_{i < \sigma} a_{\alpha,i}$ such
that $|a_{\alpha,i}| \le \lambda_i$.  Without loss of generality $\delta \in
S \Rightarrow \delta$ divisible by $\lambda^\omega$ (ordinal exponentiation).
For $\delta \in S$ let $\langle \beta^\delta_i:i < \sigma \rangle$ be
increasingly continuous with limit $\delta,\beta^\delta_i$ divisible by
$\lambda$ and  $> 0$.  For $\delta \in S$ let $\langle b^\delta_i:i < \sigma
\rangle$ be such that: $b^\delta_i \subseteq \beta^\delta_i,|b^\delta_i| \le
\lambda_i,b^\delta_i$ is increasingly continuous and $\delta = 
\dbcu_{i < \sigma} b^\delta_i$ (e.g. $b^\delta_i = \dbcu_{j_1,j_2 < i}
a_{\beta^\delta_{j_1,j_2}} \cup \lambda_i)$.  We further demand
$\beta^\delta_i \ge \lambda \Rightarrow \lambda_i \subseteq b^\delta_i 
\cap \lambda$.  Let $\langle
f^*_\alpha:\alpha < \lambda^+ \rangle$ list the two-place functions with
domain an ordinal $< \lambda^+$ and range $\subseteq \lambda^+$.  Let
$S = \dbcu_{\mu < \lambda} S_\mu$, with each $S_\mu$ stationary and
$\langle S_\mu:\mu < \lambda \rangle$ pairwise disjoint.  We now fix
$\mu < \lambda$ and will choose $\bar \gamma^\delta = \langle \gamma^\delta_i:
i < \sigma \rangle$ for $\delta \in S_\mu$ such that clause $(\alpha)$ holds
and clause $(\beta)$ holds (for every $f:\lambda^+ \rightarrow \mu$), this
clearly suffices. \nl
Now for $\delta \in S_\mu$ and $i < j < \sigma$ we can choose
$\zeta^\delta_{i,j,\varepsilon}$ (for $\varepsilon < \lambda_j)$ 
(really here we use just $\varepsilon = 0,1$) such that:
\mr
\item "{$(A)$}"  $\langle \zeta^\delta_{i,j,\varepsilon}:\varepsilon < 
\lambda_j \rangle$ is a strictly increasing sequence of ordinals
\sn
\item "{$(B)$}"  $\beta^\delta_i < \zeta^\delta_{i,j,\varepsilon} <
\beta^\delta_{i+1}$, (can even demand
$\zeta^\delta_{i,j,\varepsilon} < \beta^\delta_i + \lambda$)
\sn
\item "{$(C)$}"  $\zeta^\delta_{i,j,\varepsilon} \notin \{\zeta^\delta
_{i_1,j_1,\varepsilon}:j_1 < j,\varepsilon < \lambda_{j_1}$ (and $i_1 <
\theta$, really only $i_1 = i$ matters)$\}$
\sn
\item "{$(D)$}"  for every $\alpha_1,\alpha_2 \in b^\delta_j$, the sequence 
$\langle \text{Min}\{\lambda_j,f^*_{\alpha_1}(\alpha_2,
\zeta^\delta_{i,j,\varepsilon}):\varepsilon < \lambda_j\}\rangle$ is constant
i.e.: either
\endroster

$$
\varepsilon < \lambda_j \Rightarrow 
(\alpha_2,\zeta^\delta_{i,j,\varepsilon}) \notin \text{ Dom}(f^*_{\alpha_1})
$$

$$
\text{\ub{or} } \dsize \bigwedge_{\varepsilon < \lambda_j}
f^*_{\alpha_1}(\alpha_2,\zeta^\delta_{i,j,\varepsilon}) =
f^*_{\alpha_1}(\alpha_2,\zeta^\delta_{i,j,0})
$$

$$
\text{\ub{or} } \dsize \bigwedge_{\varepsilon < \lambda_j} f^*_{\alpha_1}
(\alpha_2,\zeta^\delta_{i,j,\varepsilon}) \ge \lambda_j.
$$
\mn
We can add $\langle f^*_{\alpha_1}(\alpha_2,\zeta^\delta_{i,j,\varepsilon}):
\varepsilon < \lambda_j \rangle$ is constant or strictly increasing (or all
undefined). \nl
Let $G = \{g:g$ a function from $\sigma$ to $2^\sigma$ such that
$(\forall i < \sigma)(i< g(i)\}$. \nl
For each function $g \in G$ we try $\bar \gamma^{g,\delta} = \langle
\zeta^\delta_{i,g(i),0},\zeta^\delta_{i,g(i),1}:i < \sigma \rangle$ i.e.
$\langle \gamma^{g,\delta}_{2i},\gamma^{g,\delta}_{2i+1} \rangle = \langle
\gamma^\delta_{i,g(i),0},\gamma^\delta_{i,g(i),1} \rangle$. \nl
Now we ask for each $g \in G$:
\enddemo 
\bn
\ub{Question}$^\mu_g$:  Does $\langle \bar \gamma^{g,\delta}:\delta \in
S_\mu \rangle$ satisfy \nl
$(\forall f \in {}^{\lambda^+}\mu)(\exists^{\text{stat}}\delta \in S_\mu)
(\dsize \bigwedge_{i < \theta}f(\gamma^{g,\delta}_{2i}) = 
f(\gamma^{g,\delta}_{2i+1}))$?.
\sn
If for some $g \in G$ the answer is yes, we are done.  Assume not, so for
each $g \in G$ we have $f_g:\lambda^+ \rightarrow \mu$ and a club $E_g$ of
$\lambda^+$ such that:

$$
\delta \in S_\mu \cap E_g \Rightarrow (\exists i < \sigma)
(f_g(\gamma^{g,\delta}_{2i}) \ne f_g(\gamma^{g,\delta}_{2i+1}))
$$
\mn
which means

$$
\delta \in S_\mu \cap E_g \Rightarrow (\exists i < \sigma)[f_g
(\zeta^\delta_{i,g(i),0}) \ne f_g(\zeta^\delta_{i,g(i),1})].
$$
\mn
Let $G = \{g_\varepsilon:\varepsilon < 2^\sigma\}$, so we can find a
2-place function $f^*$ from $\lambda^+$ to $\mu$ satisfying 
$f^*(\varepsilon,\alpha) = f_{g_\varepsilon}(\alpha)$ when $\varepsilon <
2^\sigma,\alpha < \lambda^+$.  Hence for each $\alpha < \lambda^+$ there is
$\gamma[\alpha] < \lambda^+$ such that $f^* \restriction \alpha = 
f^*_{\gamma_\alpha}$.
\mn
Let $E^* = \dbca_{\varepsilon < 2^\sigma} Eg_\varepsilon \cap \{\delta <
\lambda^+:\text{for every } \alpha < \delta \text{ we have } \gamma[\alpha]
< \delta\}$.  Clearly it is a club of $\lambda^+$, hence we can find
$\delta \in S_\mu \cap E^*$.  Now $\beta^\delta_{i+1} < \delta$ hence
$\gamma[\beta^\delta_{i+1}] < \delta$ \nl
(as $\delta \in E^*$) but $\delta =
\dbcu_{i < \sigma} b^\delta_i$ hence for some $j < \sigma,\gamma
[\beta^\delta_{i+1}] \in b^\delta_j$; as $b^\delta_j$ increases with $j$ we
can define a function $h:\sigma \rightarrow \sigma$ by $h(i) = 
\text{ Min}\{j:j > i+1$ and $\mu < \lambda_j$ and 
$\gamma[\beta^\delta_{i+1}] \in b^\delta_j\}$.
So $h \in G$ hence for some $\varepsilon(*) < 2^\theta$ we have $h =
g_{\varepsilon(*)}$.  Now looking at the choice of $\zeta^\delta_{i,h(i),0},
\zeta^\delta_{i,h(i),1}$ we know (remember $2^\sigma < \lambda_0 \le
\lambda_j$)

$$
(\forall \varepsilon < 2^\sigma)(\forall \alpha \in b^\delta_{h(i)})
(f^*_\alpha(\varepsilon,\zeta^\delta_{i,h(i),0}) = f^*_\alpha(\varepsilon,
\zeta^\delta_{i,h(i),1})).
$$
\mn
In particular this holds for $\varepsilon = \varepsilon(*),\alpha = \gamma
[\beta^\delta_{i+1}]$, so we get

$$  
f^*_{\gamma[\beta^\delta_{i+1}]}(\varepsilon(*),\zeta^\delta_{i,h(i),0}) =
f^*_{\gamma[\beta^\delta_{i+1}]}(\varepsilon(*),\zeta^\delta_{i,h(i),1}).
$$
\mn
By the choice of $f^*$ and of $\gamma[\beta^\delta_{i+1}]$ this means

$$
f_{g_{\varepsilon(*)}}(\zeta^\delta_{i,h(i),0}) =
f_{g_{\varepsilon(*)}}(\zeta^\delta_{i,h(i),1)})
$$
\mn
but 
$h = g_{\varepsilon(*)}$ and the above equality means 
$f^*_{g_{\varepsilon(*)}}
(\gamma^{g_{\varepsilon(*)},\delta}_{2i}) = f^*_{g_{\varepsilon(*)}}
(\gamma^{g_{\varepsilon(*)},\delta}_{2i+1})$, and this holds for every
$i < \sigma$, and $\delta \in E^* \Rightarrow \delta \in 
E_{g_{\varepsilon(*)}}$ so we get a contradiction to the choice of
$(f_{g_{\varepsilon(*)}},E_{\varepsilon(*)})$. \nl
So we have finished proving $(\alpha) + (\beta)$.
\mn
How do we get $(\beta)^+$ of \scite{1.1}, too? \nl

The first difference is in phrasing the question, now it is, for $g \in G$:
\mn
Question$^\mu_g$:  Does $\langle \bar \gamma^{g,\delta}:\delta \in
S_\mu \rangle$ satisfies:

$$
\biggl( (\forall f_0 \in {}^{\lambda^+}\mu_0)(\forall f_1 \in {}^{\lambda^+}
\mu_1) \ldots (\forall f_i \in {}^{\lambda^+} \mu_i) \ldots 
\biggr)_{i < \sigma}
(\exists^{\text{stat}}\delta \in S_\mu) (\dsize \bigwedge_{i < \sigma}
f_i(\gamma^{g,\delta}_{2_i}) = f_i(\gamma^{g,\delta}_{2i+1})).
$$
\mn
If for some $g$ the answer is yes, we are done, so assume not so we have
$f_{g,i} \in {}^{\lambda^+}(\mu_i)$ for $g \in G,i < \sigma$ and club $E_g$
of $\lambda^+$ such that

$$
\delta \in S_\mu \cap E_g \Rightarrow (\exists i < \sigma)(f_{g,i}
(\gamma^{g,\delta}_{2i}) \ne f_{g,i}(\gamma^{g,\delta}_{2i+1})).
$$
\mn
A second difference is the choice of $f^*$ as $f^*(\sigma \varepsilon +i,
\alpha) = f_{g_\varepsilon,i}(\alpha)$ for $\varepsilon < 2^\sigma$, \nl
$i < \sigma,\alpha < \lambda^+$. \nl
Lastly, the equations later change slightly.
\hfill$\square_{\scite{1.1}}$
\bn
\ub{\stag{1.2} Fact}:  Under the assumptions of \scite{1.1} letting
$\bar \lambda = \langle \lambda_i:i < \sigma \rangle$ be increasingly
continuous with limit $\lambda$ we have
\mr
\item "{$(*)_1$}"  we can find $\langle <\gamma^\delta_\zeta:
\zeta < \lambda>:\delta \in S \rangle$ such that
{\roster
\itemitem{ $(\alpha)$ }  $\gamma^\delta_i$ is increasing in $i$ with limit
$\delta$
\sn
\itemitem{ $(\beta)^+$ }  if $f_i:\lambda^+ \rightarrow
\lambda_{i+1}$, for $i < \sigma$, \ub{then} the following set is stationary
$\{\delta \in S:f_i(\gamma^\delta_\zeta) = f_i(\gamma^\delta_\xi)$ when
$\zeta,\xi \in [\lambda_i,\lambda_{i+1})\}$
\endroster}
\item "{$(*)_2$}"  moreover if $F_i:[\lambda^+]^{< \lambda} \rightarrow
[\lambda^+]^{\lambda^+}$ for $i < \sigma$, \ub{then} we can demand
{\roster
\itemitem{ $(i)$ }  $\{\gamma^\delta_\zeta:\zeta \in [\lambda_i,\lambda_{i+1}]
\} \subseteq F_i(\{\gamma^\delta_\zeta:\zeta < \lambda_i\})$,
\sn
\itemitem{ $(ii)$ }  $|\{\langle \gamma^\delta_\zeta:\zeta < \zeta^* \rangle:
\gamma^\delta_{\zeta^*} = \gamma\}| \le \lambda$ for each $\gamma < \lambda^+$
and $\zeta^* < \sigma$
\endroster}
\item "{$(*)_3$}"  moreover, if $\langle C_\delta:\delta \in S \rangle$ is
given, it guess clubs, $C_\delta = \{\alpha[\delta,i]:i < \sigma\},
\alpha[\delta,i]$ divisible by $\lambda^\omega$ increasing in $i$ with 
limit $\delta,\langle \text{cf}(\alpha[\delta,i+1]):i < \sigma \rangle$ is
increasing with limit $\lambda$ letting $\beta(\delta,i) =
\dsize \sum_{j<i} \lambda_j \times \text{ cf}(\alpha[\delta,j])$ we can demand
sup$\{\gamma^\delta_i:\gamma^\delta_i < \beta[\delta,j]\} = \alpha[\delta,j]$
(and in $(\beta)^+$ the division to intervals is by the $\beta[\delta,j]$
not $[\lambda_i,\lambda_{i+1}]$) and for every $f_i \in {}^{(\lambda^+)}
(\mu_i)$ for $i < \sigma$ where $\mu_i < \lambda$ and club $E$ of $\lambda^+$,
for stationarily many $\delta \in S$ we have $\{\gamma^\delta_i:i < \lambda\}
\subseteq E$ and $f_i(\gamma^\delta_\zeta) = f_i(\gamma^\delta_\varepsilon)$,
when $\zeta,\varepsilon \in [\beta[\delta,i],\beta[\delta,i+1])$.
\endroster
\bigskip

\demo{Proof}  The same proof as in \scite{1.1} for $(*)_1$.

Minor changes for $(*)_2,(*)_3$.  For clause (ii) of $(*)_2$ note that we
can find \nl
$\langle {\Cal P}_\alpha:\alpha < \lambda^+ \rangle$ such that:
\mr
\item "{$(a)$}"  ${\Cal P}_\alpha \subseteq [\alpha]^{< \lambda}$,
\sn
\item "{$(b)$}"  if $A \subseteq \alpha$ then for some $\bar A = \langle
A_i:i < \sigma \rangle$ we have: \nl
$A_i$ is increasingly continuous, $A_i \in {\Cal P}_\alpha$ and $A =
\dbcu_{i < \sigma} A_i$
\sn
\item "{$(c)$}"  ${\Cal P}_\alpha$ is closed under union of $< \sigma$
members.
\endroster  
\enddemo  
\newpage

\head {\S2 Case C:Forcing for successor of singular} \endhead  \resetall
\bn
See \cite{Sh:587}(xxx).
\mn
\demo{\stag{2.1} Hypothesis}  1) $\lambda$ strong limit singular
$\sigma = \text{ cf}(\lambda) < \lambda,\kappa = \lambda^+,\mu^* \ge \kappa,
2^\lambda = \lambda^+$. \nl
2) $\hat{\Cal E}_0 \subseteq \{\bar a:\bar a = \langle a_i:i \le \delta
\rangle,a_i \in [\mu^*]^{< \kappa}$ increasingly
continuous, $\lambda \subseteq a_0\}$ non-trivial (i.e. for $\theta =
\text{ cf}(\theta) < \lambda,\chi$ large enough and $x \in {\Cal H}(\chi)$ 
we can find
$\langle N_i:i \le \theta \rangle$ obeying $\bar a \in \hat{\Cal E}_0$ (with
error some $n$ see \cite[xx]{Sh:587} and such that $x \in N_0$). \nl
3) $\hat{\Cal E}_1 \subseteq \{\bar a:\bar a = \langle a_i:i \le \sigma
\rangle,a_i$ increasingly continuous $|a_i| < \lambda,\lambda +1 \subseteq
\dbcu_{i < \sigma} a_i\}$.
\enddemo
\bigskip

\definition{\stag{2.2} Definition}  $(\hat{\Cal E}_0,\hat{\Cal E}_1)$
non-trivial if:
\mr
\item "{$(a)$}"  $\hat{\Cal E}_0$ is non-trivial
\sn
\item "{$(b)$}"  for $\chi$ large enough and 
$x \in {\Cal H}(\chi)$, we can find $\bar M = \langle M_i:
i \le \sigma \rangle$ and \nl
$\langle \bar N^i:i < \sigma \rangle$ such that
{\roster
\itemitem{ $(\alpha)$ }  $M_i \prec ({\Cal H}(\chi),\in,<^*_\chi)$ and
\sn
\itemitem{ $(\beta)$ }  $\bar M$ obeys some $\bar a \in \hat{\Cal E}_1$ for
some finite error and
\sn
\itemitem{ $(\gamma)$ }   $[M_{i+1}]^{2^{\|M_i\|}} \subseteq M_{i+1}$
\sn
\itemitem{ $(\delta)$ }  $N^i_\alpha \prec ({\Cal H}(\chi),\in,<^*_\chi)$
and
\sn
\itemitem{ $(\varepsilon)$ }  $\langle N^i_\alpha:\alpha \le \delta_i 
\rangle = \bar N^i$ obeys some $\bar b_i \in \hat{\Cal E}_0$ for some finite
error 
\sn
\itemitem{ $(\zeta)$ }  $M_{i+1} \subseteq N^i_{\delta_i}$
\sn
\itemitem{ $(\eta)$ }  $\delta_i \subseteq M_{i+1}$,
\sn
\itemitem{ $(\theta)$ }  cf$(\delta_i) > 2^{\|M_i\|}$, 
\sn
\itemitem{ $(\iota)$ }  $\alpha < \delta_i \Rightarrow 
\bar N^i \restriction \alpha \in M_{i+1}$
\sn
\itemitem{ $(\kappa)$ }  $N^i_{\delta_i} \prec N^j_0$ for $i<j$
\sn
\itemitem{ $(\lambda)$ }   $M_i \prec N^i_0,M_i \in N^i_0$
\sn
\itemitem{ $(\mu)$ }  if $\beta < \sigma$ is a limit ordinal, \ub{then}
$\langle N^i_{\delta_i}:i < \beta \rangle \char 94 \langle 
\dbcu_{i < \beta} N^i_{\delta_i} \rangle$ obeys some 
$\bar a^\beta \in \hat{\Cal E}_0$ (needed for the ghost coordinates).
\endroster}
\endroster
\enddefinition
\bigskip

\remark{Remark}  1) Letting $\bar N = \bar N^0 \char 94 \bar N^1 \ldots$ and
$N_\varepsilon = N^i_\alpha$ if $\varepsilon = \dsize \sum_{j<i} \delta_j
+ \alpha$ then $\ell g(\bar N) = \lambda$ and $\bar N \restriction (i_0+1)
\in N_{i_0+1}$ so $\bar N$ is $\prec$-increasingly continuous, $\bar N
\restriction \gamma \in N_{i+1}$.  \ub{Note} $\dbcu_i N^i_{\delta_i} =
\dbcu_i N^i_0 = M_\sigma$. \nl
2) [\ub{Variant}:  
$\delta_i = \sup(\delta_i \cap M_{i+1}),\alpha \in \delta_i
\cap M_{i+1} \Rightarrow \bar N^i \restriction \alpha \in M_{i+1}$.
Saharon check.] \nl
3) Clearly $(b) \Rightarrow (a)$ in Definition \scite{2.2}.
\endremark
\bigskip

\definition{\stag{2.3} Definition}  $Q$ is [strongly] 
$\hat{\Cal E}_0$-complete: as in the inaccessible case.
\enddefinition
\bigskip

\proclaim{\stag{2.4} Claim}  1) Assume $\hat{\Cal E}_0$ is non-trivial
and we have $(< \kappa)$-support iteration $\bar Q$, and
(non-trivial) $\hat{\Cal E}_0$, $\Vdash_{P_i} ``{\underset\tilde {}\to Q_i}$ 
is $\hat{\Cal E}_0$-complete". \nl
\ub{Then} $P_\gamma$ is $\hat{\Cal E}_0$-complete 
(hence $P_\gamma/P_\beta$). \nl
2) $Q$ is $\hat{\Cal E}_0$-complete $\Rightarrow V^Q \models \hat{\Cal E}_0$
non-trivial.
\endproclaim
\bigskip

\demo{Proof}  Just like the inaccessible (here the choice ``for any regular
cardinal $\theta < \kappa$" rather than ``for any cardinal $\theta <
\kappa$" is important).
\enddemo
\bigskip

\definition{\stag{2.5} Definition}  $Q$ is $(\hat{\Cal E}_0,
\hat{\Cal E}_1)$-complete if: 
\mr
\item "{$(A)$}"  $Q$ is strongly ${\Cal E}_0$-complete
\sn
\item "{$(B)$}"  For some $x$ if $x \in {\Cal H}(\chi)$ and 
$\bar M,\langle \bar N^i:i < \sigma \rangle$ are as in Definition 
\scite{2.2}, then ($Q \in M_0$, and)
{\roster
\itemitem{ $(*)$ }  for any $i_0 < \sigma$ and $r \in M_{i_0+1} \cap Q$, in 
the following game between the players COM and INC, the player INC has no
winning strategy, we index the moves by $i \in [i_0,\sigma)$, 
\ub{in the ith move} the player COM chooses $p_i \in Q \cap M_{i+1}$
such that $r \le p_i$ and for $j \in [i_0,i),p_i$ is an upper bound of
$\bar q^j$.  Then the INC player chooses $\alpha_i < \delta_i$ and $\bar q' =
\langle q^i_\alpha:\alpha < \delta_i,\alpha \ge \alpha_i \rangle$ such that
$p_i \le p^i_{\alpha_i}$ and $\bar q^i$ is generic for $\bar N^i$ and
$\dsize \bigwedge_{\alpha \in \delta_i \cap M_{i+1}} \bar q^i \restriction
\alpha \in M_{i+1}$. \nl
We may phrase it as one game, where INC chooses $i_0 < \delta$ and $r \in
M_{i_0+1} \cap Q$.
\endroster}
\endroster
\enddefinition
\bigskip

\proclaim{\stag{2.6} Claim}  Assume $\hat{\Cal E}_0$ is non-trivial
(see \scite{2.1}(2)).  If ${\Cal T}$ is a $(w,\alpha)$-tree,
$w \subseteq r = \ell g(\bar Q),|w| < \lambda,\bar Q$ is $(< \kappa)$-support
iteration of $\hat{\Cal E}_0$-complete forcing notions and 
$\langle p_t:t \in {\Cal T} \rangle \in TFr(\bar Q),\langle
{\Cal I}_\alpha:\alpha \in w \cup \{\gamma\} \rangle,{\Cal I}_\alpha \subseteq
P_\alpha$ dense, \ub{then} there is $\langle q_t:t \in {\Cal T} \rangle \in
TFr(\bar Q)$ above $\langle p_t:t \in {\Cal T} \rangle$ such that
$t \in {\Cal T} \Rightarrow q_t \in {\Cal I}_{\text{rk}(t)}$.
\endproclaim
\bigskip

\demo{Proof}  The point is $|{\Cal T}| < \lambda$, otherwise as in 
inaccessible case in \cite{Sh:587}.
\enddemo
\bigskip

\proclaim{\stag{2.7} Claim}  Assume $(\hat{\Cal E}_0,\hat{\Cal E}_1)$ is
non-trivial. \nl
1) If $\bar Q$ is $(< \kappa)$-support
$\Vdash_{P_i} ``Q_i$ is $(\hat{\Cal E}_0,\hat{\Cal E}_1)$-complete",
\ub{then} $P_\gamma$ is $(\hat{\Cal E}_0,\hat{\Cal E}_1)$-complete
$(\hat{\Cal E}_0,\hat{\Cal E}_1)$ non-trivial. \nl
2) If $Q$ is $(\hat{\Cal E}_0,\hat{\Cal E}_1)$-complete, then also in $V^Q$.
\endproclaim
\bigskip

\demo{Proof}  1) Like in the inaccessible case (\cite{Sh:587}) with 
$\bar M,\langle \bar N^i:i < \sigma \rangle$ as in Definition \scite{2.2}.  
We define the trees
point: in stage i using trees ${\Cal T}_i$ with set of levels $w_i = M_i
\cap \gamma$ and looking at all possible moves of COM, i.e. $p_i \in
M_{i+1} \cap P_\gamma$, so constructing this tree of conditions in
$\delta_i$ stages, in stage $\varepsilon < \delta_i$, has 
$|N^i_\varepsilon \cap M_{i+1} \cap Q|^{2^{\|M_i\|}}$ nodes.
\sn
Now
\mr
\item "{{}}"  $p \in P_\gamma \cap M_{i+1} \nRightarrow \text{ Dom }
p \subseteq M_{i+1}$ but
\sn
\item "{{}}"  $p \in P_\gamma \cap M_{i+1} \Rightarrow \text{ Dom } p
\subseteq M_\sigma = \dbcu_{i < \sigma} N^i_{\delta_i}$ \nl
$p \in P_\gamma \cap N^i_\varepsilon \Rightarrow \text{ Dom } p \subseteq
N^i_{\varepsilon + 1}$.
\ermn
So in limit cases $i < \sigma$: the existence of limit is by the clause
$(\mu)$ of Definition \scite{2.2}.
In the end we use the winning of the play
and then need to find a branch in the tree of conditions of level $\sigma$:
like Case A using $\hat{\Cal E}_0$.  \hfill$\square_{\scite{2.7}}$
\enddemo
\bn
\ub{\stag{2.8} Generalization}  $\hat{\Cal E}_1$ is a set of triples
$\langle \bar a,\langle \bar b^i:i < \sigma \rangle,\bar \lambda \rangle,
\bar a = \langle a_i:i \le \sigma \rangle$ as before
$\bar b^i = \langle b^i_\alpha:\alpha \le \delta_i \rangle \in
\hat{\Cal E}_0,b^i_{\delta_i} = a_{i+1},a_i \subseteq b^i_0,\lambda = 
\langle \lambda_i:i < \sigma \rangle$ an increasing sequence of cardinals
$< \lambda,\sum \lambda_i = \lambda$ and $(\bar M,\langle \bar N^i:i < \sigma
\rangle$ obeys $(\bar a,\langle \bar b^i: i < \bar \lambda \rangle$ is
$M_i \cap \mu^* = a_i,\bar N^i$ obeys $\bar b^i$ all things in \scite{2.2}
but $\lambda_i \ge \|M_i\|,\lambda_i \ge \dsize \prod_{j \le i} \|M_j\|,
[N^i_\alpha]^{\lambda_i} \subseteq N^i_{\alpha +1}$ for $\alpha < \delta_i$
(so earlier $\lambda_i = 2^{\|M_i\|}$).
\bn
\ub{\stag{2.9} Conclusion}  Assume
\mr
\item "{$(a)$}"  $S \subseteq \{\delta < \kappa:\text{cf}(\delta) =
\sigma\}$ is stationary not reflecting
\sn
\item "{$(b)$}"  $\bar a = \langle \bar a_\delta:\delta \in S \rangle,
\bar a_\delta = \langle a_{\delta,i}:i \le \sigma \rangle,\delta = \cup
a_{\delta,i}$ \nl
[variant note: $\bar \lambda^\delta = \langle \lambda^\delta_i:i < \sigma
\rangle$ increasing with limit $\lambda$]
\sn
\item "{$(c)$}" letting $\mu^* = \kappa,\hat{\Cal E}_0 = \hat{\Cal E}_0[S] =
\{\bar a:a_i \in \kappa \backslash S$ increasing continuous$\}$ \nl
$\hat{\Cal E}_1 = \{\bar a_\delta:\delta \in S\}$ \nl
(or $\{\langle \bar a_\delta,\langle \bar b^{i,\delta}:i < \delta \rangle,
\bar \lambda^\delta \rangle:\delta \in S\}$ appropriate for (\scite{2.8}).
\nl
We assume the pair $(\hat{\Cal E}_0,\hat{\Cal E}_1)$ is not trivial
\sn
\item "{$(d)$}"  $\mu = \mu^\kappa,\kappa < \tau = \text{ cf}(\tau) < \mu$.
\ermn
\ub{Then} for some $(\hat{\Cal E}_0,\hat{\Cal E}_1)$-complete $P$ of
cardinality $\mu$

$$
\align
\Vdash_P ``&\text{forcing axioms for } (\hat{\Cal E}_0,\hat{\Cal E}_1)
\text{-complete forcing notion} \\
  &\text{of cardinality } \le \kappa \text{ and } < \tau
\text{ of open dense sets}" \text{ and } S \text{ is still} \\
  &\text{stationary (by preservation of } (\hat{\Cal E}_0,\hat{\Cal E}_1)
\text{-nontrivial}.
\endalign
$$
\bigskip

\demo{Proof}  As in the inaccessible case.

In particular the $\kappa^+$-c.c. lemma work.
\enddemo   
\bn
\ub{\stag{2.10} Application}: In $V^P$ of \scite{2.9}:
\mr
\item "{$(a)$}"  if $\theta < \lambda,A_\delta \subseteq \delta = \sup
(A_\delta)$ for $\delta \in S,|A_\delta| \le \theta,(\forall^* i < \sigma)
[a_\delta \cap \sup a_{\delta,i} \subseteq a_{\delta,i}],\bar h =
\langle h_\delta:\delta \in S \rangle,h_\delta:A \rightarrow \theta$,
\ub{then} for some $h:\kappa \rightarrow \theta$, for some club $E$ of
$\kappa,(\forall \delta \in S \cap E)[h_\delta \subseteq^* h]$ bounded set
of errors where $h' \subseteq^* h''$ means sup(Dom $h') > \sup \{\alpha_i,
\alpha \in \text{ Dom } h'$ and $\alpha \notin \text{ Dom}(h'')$ or \nl
$\alpha \in \text{ Dom}(h'') \and h'(\alpha) \ne h''(\alpha)\}$
\sn
\item "{$(b)$}"  if we add: ``$h_\delta$ constant", \ub{then} we can omit
the assumption \nl
$(\forall^* i < \sigma)(A_\delta \cap \sup a_{\delta,i} 
\subseteq a_{\delta,i})$
\sn
\item "{$(c)$}"  we can weaken $|A_\delta| < \theta$ to $|A_\delta \cap
a_{\delta,i+1}| \le |a_{\delta,i}|$
\sn
\item "{$(d)$}"  in (c) we can weaken $|A_\delta| < \theta \vee
|A_\delta \cap a_{\delta,i+1}| \le |a_{\delta,i}|$ to 
$h_\delta \restriction a_{\delta,i+1}$ belongs to $M_{i+1} \cap N^i_\alpha$
for some $\alpha < \delta_i$ \nl
(remember cf(sup $a_{\delta,i+1}) > \lambda^\delta_i$).
\endroster
\bigskip

\remark{\stag{2.10B} Remark}  1) Compared to \cite{Sh:186} the new 
point is (b). \nl
2) You may complain why not having the best of (a) + (b), i.e. combine their
good points.  The reason is that this is impossible by \S1, \S4 (c), 
(d) - similar situation in the inaccessible.
\endremark
\bigskip

\demo{Proof}  Should be clear.  Still we say something in case $h_\delta$ 
constant (b). \nl
Let

$$
\align
Q = \{(h,C):&\alpha =: \text{ Dom}(h) < \kappa = \lambda^+, \\
  &C \text{ a closed subset of } \alpha +1,\alpha \in C \\
  &(\forall \delta \in C \cap S)(h_\delta \subseteq^* h)\}.
\endalign
$$
\enddemo
\bn
So \nl
\ub{Fact}:  $(h,C) \in Q,\alpha = \text{ Dom}(h) < \beta \in \kappa
\backslash S,\gamma < \theta \Rightarrow (h,C) \le (h \cup 
\gamma_{[\alpha,\beta]},C \cup \{\beta\}) \in Q$. \nl
In the game from \scite{2.2}, we can even prove that the player COM has a
winning strategy: in stage i (non-trivial): if $h_\delta$ is constantly
$\gamma < \theta$ or just $h_\delta \restriction (A_\delta \cap 
a_{\delta,i+1} \backslash a_{\delta,i})$ is constantly $\gamma < \theta$
then we let

$$
\align
p_i = \biggl( &\dbcu_{\Sb j < i \\ \zeta < \delta_i \endSb}
\text{ the function } q^j_\zeta) \cup 
\gamma_{[\text{Dom } \dbcu_{\Sb j < i \\ \zeta < \delta \endSb} q^j_\zeta,
\beta_i)}, \\
  &\text{closure} \bigl(\dbcu_{\Sb j < i \\ \zeta < \delta_i \endSb}
(\text{club of } q^j_\zeta) \cup \{\beta\} \bigl) \biggl)
\endalign
$$
\mn
for $\beta \in M_{i+1} \cap \kappa \backslash M_i \text{ large enough such
that } A_\delta \cap M_{i+1} \cap \kappa \subseteq \beta$.
\bigskip

\remark{Remark}  In the example of uniformizing (see \cite{Sh:587}) 
if we use this forcing, the density is less problematic.
\endremark
\bigskip

\proclaim{\stag{2.11} Claim}  In \scite{2.10}'s conclusion we can omit the
club $E$ and demand \nl
$(\forall \delta \in S)(h_\delta \subseteq^* h)$ provided
that we add in \scite{2.9} (recalling $S \subseteq \kappa$ does not reflect)
is a set of limits and

$$
\bar A = \langle A_\delta:\delta \in S \rangle,A_\delta \subseteq \delta =
\sup(A_\delta)
$$
\mn
satisfies
\mr
\item "{$(*)$}" $\delta_1 \ne \delta_2 \text{ in } S \Rightarrow \sup
(A_{\delta_i} \cap A_\delta) < \delta_1 \cap \delta_2$.
\ermn
(If $(\forall \delta \in S)(\text{cf}(\delta) = \theta) \and
\dsize \bigwedge_\delta \text{ otp}(A_\delta) = \theta$ this always holds.
\endproclaim
\bigskip

\demo{Proof}  We define $Q = \{h:\alpha = \text{ Dom}(h) < \kappa$, Rang
$h \subseteq \dbcu_{\delta \in S} \text{ Rang } h_\delta\}$ ordered by
$\subseteq$.  Now we can have the parallel of the fact that: 
if $p:\alpha \rightarrow$ (set of colours), $\alpha < \beta < \kappa$. \nl
We can find $\langle A'_\delta:\delta \in S \cap (\alpha +1) \rangle$ such
that $A'_\delta \subseteq A_\delta,\sup(A^*_\delta \backslash A'_\delta)
< \delta$, \nl
$\bar A' = \langle A'_\delta:\delta \in S \cap (\alpha + 1) \rangle$
pairwise disjoint.
\sn
Now choose $q$ as follows

$$
\text{Dom}(q) = \beta
$$

$$
q(j) = \cases p(j) &\text{ \ub{if} } \quad j < \alpha \\
h_\delta(j) &\text{ \ub{if} } \quad j \in A'_\delta \backslash \alpha,
\delta \in S \cap (\beta +1) \backslash (\alpha +1) \\
0 &\text{ \ub{if} otherwise}. \endcases
$$
\mn
Why does $\bar A'$ exist?  Prove by induction on $\beta$ that any
$\langle A'_\delta:\delta \in S \cap (\alpha +1) \rangle,\alpha < \beta$,
can we end extending to $\langle A'_\delta:\delta \in S \cap (\beta +1)
\rangle$.
\bigskip

\remark{\stag{2.12} Remark}  Note: concerning $\kappa$ inaccessible we could
immitate what is here: having $M_{i+1} 
\underset\ne {}\to \prec N^i_{\delta_i},
\dbcu_{i < \delta} M_i = \dbcu_{i < \delta} N^i_{\delta_i}$.

As long as we are looking for a proof no sequence of length $< \kappa$ are
added, the gain is meagre (restricting the $\bar q$'s by $\bar q
\restriction \alpha \in N'_{\alpha +1}$).  Still if you want to make the
uniformization and some diamond may consider.
\endremark
\bn
\ub{\stag{2.13} Comment}:  We can weaken further the demand, by letting COM
have more influence.  E.g. we have (in \scite{2.2}) $\delta_i = \lambda_i
= \text{ cf}(\lambda_i) = \|M_{i+1}\|,D_i$ a $|a_i|^+$-complete filter on
$\lambda_i$, the choice of $\bar q^i$ in the result of a game in which INC
should have chose a set of player $\in D_i$ and $\diamondsuit_{D_i}$ holds
(as in the treatment of case $E^*$ here).  More details?
\sn
The changes are obvious, but I do not see an application at the moment.
\bigskip

\proclaim{\stag{2.14} Claim}  In \scite{2.10}(9) we get:
\mr
\item "{$(*)$}"  if for $\delta \in S,A_\delta \subseteq \delta,|A_\delta|
< \lambda$ then for some club $E$ of $\kappa$ we have \nl
$\delta \in E \cap S \Rightarrow \text{ otp}(A_\delta \cap E) \le \sigma$.
\endroster
\endproclaim 
\bigskip

\demo{Proof}  Straight.
\enddemo
\enddemo
\newpage

\head {\S3 $\kappa^+$-c.c. and $\kappa^+$-pic} \endhead  \resetall
\bn
We intend to generalize pic of \cite[Ch.VIII,\S1]{Sh:f}. The intended use
is for iteration
with each forcing $> \kappa$ - see use in \cite{Sh:f}.  (At present we need
each $Q_i$ of cardinality $\le \kappa$.  Usually $\mu = \kappa^+$.)
\definition{\stag{3.1} Definition}  Assume:
\mr
\item "{$(a)$}"   $\mu = \text{ cf}(\mu) > |\alpha|^{< \kappa}$ for 
$\alpha < \mu$
\sn
\item "{$(b)$}"  $(\kappa,\mu^*,\hat{\Cal E}_0)$ as usual in \cite{Sh:587}
that is $\kappa = \text{ cf}(\kappa) > \aleph_0,\mu^* \ge \kappa,
\hat{\Cal E}_0 \subseteq \{\bar a:\bar a$ an increasing continuous sequence of
members of $[\mu^*]^{< \kappa}$ of limit length $< \kappa\}$ and
\sn
\item "{$(c)$}"  $S^\square \subseteq 
\{\delta < \mu:\text{cf}(\delta) \ge \kappa\}$ stationary.
\ermn
We say $Q$ satisfies $(\mu,S^\square,\hat{\Cal E}_0)$-pic if: for some 
$x \in {\Cal H}(\chi)$ (can be omitted, essentially)
\mr
\item "{$(*)$}"  if $\langle \mu,S^0,\hat{\Cal E}_0,x
\rangle \in N^\alpha_0$ for $\alpha \in S^\square,\bar N^\alpha = \langle
\bar N^\alpha_i:i \le \delta_\alpha \rangle$ obeys $\bar a^\alpha \in
\hat{\Cal E}_0$ (with some error $n$) (so here we have 
$\|N^\alpha_{\delta_\alpha}\| < \kappa,\delta_\alpha < \kappa$) and
$\bar p^\alpha$ is generic for $\bar N^\alpha$ and for every $i$, if
$\alpha \in S^\square,i < \delta_\alpha$ then $\langle (\bar N^\beta
\restriction (i+1),\bar p^\beta \restriction (i+1)):\beta \in S^\square
\rangle$ belongs to $N^\alpha_{i+1}$ and we define a function $g$ with domain
$S^\square,g(\alpha) = (g_0(\alpha),g_1(\alpha)),g_0(\alpha) = N^\alpha
_{\delta_\alpha} \cap (\dbcu_{\beta < \alpha} N^\beta_{\delta_\beta}),
g_1(\alpha) = (N^\alpha_{\delta_1},N^\alpha_i,c)_{i < \delta_1,c \in 
g_0(\alpha)}/ \cong$, \ub{then} we can find a club $C$ of $\mu$ and a 
regressive function $g:C \cap S \rightarrow \mu$ such that: \nl
if $\alpha < \beta \and g(\alpha) = g(\beta) \and \alpha \in 
C \cap S \and \beta \in C \cap S$ \ub{then} for some 
$h,N^\alpha_{\delta_\alpha}
\underset h {}\to \cong N^\beta_{h_\alpha}$ (really unique), $h$ maps
$N^\alpha_i$ to $N^\beta_i,p^\alpha_i$ to $p^\beta_i$ and
$\{p^\alpha_i:i < \delta_\alpha\} \cup \{p^\beta_i:i < \delta_\beta\}$ has
an upper bound.
\endroster
\enddefinition
\bigskip

\proclaim{\stag{3.1A} Claim}  Assume (a),(b) of \scite{3.1} and
\mr
\item "{$(c)$}"  $\hat{\Cal E}_0$ is non-trivial, that is for every $\alpha$
large enough and $x \in {\Cal H}(\chi)$ there is
$\bar N = \langle N_i:i \le \delta \rangle$ increasingly continuous,
$N_i \prec ({\Cal H}(\chi),\in),x \in N_i,\|N_i\| < \kappa,\bar N \restriction
(i+1) \in N_{i+1}$ and $\bar N$ obeys some $\bar a \in \hat{\Cal E}_0$ with
some finite error $n$)
\sn
\item "{$(d)$}"  $Q$ is an $\hat{\Cal E}_0$-complete forcing notion not 
adding bounded subsets of $\kappa$ and 
$Q$ satisfies $(\mu,S^\square,\hat{\Cal E}_0$)-pic.
\ermn
\ub{Then} $Q$ satisfies the $\mu$-c.c. (and for every $x \in {\Cal H}(\chi),
\chi$ large enough, $p \in Q$ for some $\bar N$ as in clause (c),
$\langle Q,p,s \rangle \in N_0$ and for some increasing sequence $\bar p =
\langle p_i:i \le \ell g(\bar N) \rangle$ of members of $Q,\bar p$ is generic
for $\bar N$ and $p \le p_0$).
\endproclaim
\bigskip

\demo{Proof}  Straight.
\enddemo
\bigskip

\proclaim{\stag{3.1B} Claim}  If Definition \scite{3.1}, we can allow Dom$(g)$
to be a subset of $\lambda,\langle A_i:i < \mu \rangle$ be an increasingly
continuous sequence of sets, $|A_i| < \mu,N^\alpha_{\delta_\alpha} \subseteq
A_{\alpha +1}$ replacing the definition of $g$ by $g_0(\alpha) =
N^\alpha_{\delta_\alpha} \cap A_\alpha,g_i(\alpha) = 
(N^\alpha_{\delta_\alpha},N^\alpha_i,c)_{i < \delta_\alpha,c \in g_0(c)} /
\cong$ (and get equivalent definition).
\endproclaim
\bigskip

\remark{Remark}  If Dom$(g) \cap S^\square$ is not stationary, the definition
says nothing.
\endremark
\bigskip

\demo{Proof}  Straight.
\enddemo
\bigskip

\proclaim{\stag{3.2} Claim}  Assume clauses (a), (b) of \scite{3.1} and (c)
of \scite{3.1A}.

For $(< \kappa)$-support iteration $\bar Q$ if we have 
$\Vdash_{P_i} ``{\underset\tilde {}\to Q_i}$ is $\hat{\Cal E}_0$-complete
\footnote{this includes the facts of initial segments of $\bar a \in
\hat{\Cal E}_0$}, $(\mu,S^\square,\hat{\Cal E}_0)$-pic" and
forcing with Lim$(\bar Q)$ add no bounded subsets of
$\kappa$,  \ub{then} $P_\gamma$ and $P_\gamma/P_\beta$, for $\beta \le \gamma
\le \ell g(\bar Q)$ are $\hat{\Cal E}_0$-complete
$(\mu,S^\square,\hat{\Cal E}_0)$-pic.
\endproclaim
\bigskip

\demo{Proof}  Similar to \cite[Ch.VIII]{Sh:f}.

We prove this by induction.  Let $\Vdash_{P_i} ``{\underset\tilde {}\to Q_i}$
is $(\mu,S^\square,\hat{\Cal E}_0)$-pic as witnessed by
${\underset\tilde {}\to x_i}$ and let ${\underset\tilde {}\to \chi_i} = 
\text{ Min}\{\chi:{\underset\tilde {}\to x_i} \in {\Cal H}(\chi)\}"$.

Let $x = (\mu^*,\kappa,\mu,S^\square,\hat{\Cal E}_0,\langle
({\underset\tilde {}\to \chi_i},{\underset\tilde {}\to x_i}):i < \ell g
(\bar Q) \rangle)$ and assume $\chi$ is large enough such that $x \in
{\Cal H}(\chi)$ and let $\langle (\bar N^\alpha,\bar p^\alpha):\alpha \in
S^\square \rangle$ be as in Definition \scite{3.1}, $\bar N^\alpha = \langle
N^\alpha_i:i \le \delta_\alpha \rangle$.  We define a $h,C$ such that
\mr
\item "{$\boxtimes_1$}"  $C$ is a club of $\mu,g$ is a function with domain
$S^\square$ \nl
$g(\alpha) = \langle g_\ell(\alpha):\ell < 2 \rangle$ \nl
$g_0(\alpha) = (N^\alpha_{\delta_\alpha}) \cap (\dbcu_{\beta < \alpha}
N^\beta_{\delta_\beta})$ \nl
$g_1(\alpha) =$ the isomorphic type of $(N^\beta_{\delta_\beta},N^\beta_i,
p^\beta_i,c)_{c \in h_1(\alpha)}$
\sn
\item "{$\boxtimes_2$}"  $C = \{\delta < \mu:\text{if } \alpha \in
[\delta,\mu),h_0(\alpha) \subseteq \dbcu_{\beta < \delta} N^\beta
_{\delta_\alpha}$ then for some $\alpha' < \delta,h(\alpha') = h(\alpha)\}$.
\ermn
Fix $y$ such that $S_y = \{\alpha \in S^\square:g(\alpha) = y \text{ and }
\alpha \in C\}$ is stationary.

Let $w_\alpha = \dbcu_{i < \delta_\alpha} \text{ Dom}(p^\alpha_i),
w^* = w_\alpha \cap g_0(\alpha)$ for $\alpha \in S_y$ (the set does not
depend on the $\alpha$).  For each $\zeta \in w^*_y$ we define a
$P_\zeta$-name, ${\underset\tilde {}\to S_{y,\zeta}}$ as follows:

$$
{\underset\tilde {}\to S_{y,\zeta}} = \{\alpha \in S_y:(\forall i <
\delta_\alpha)(p^\alpha_i \in {\underset\tilde {}\to G_{P_\zeta}})\}.
$$
\mn
Now apply Definition \scite{3.1} in $V^{P_\zeta}$ to
$\left < ( \langle N^\alpha_i[{\underset\tilde {}\to G_{P_\zeta}}]:i \le
\delta_\alpha \rangle,\langle p^\alpha_i(\zeta)
[{\underset\tilde {}\to G_{P_\zeta}}]:c < \delta_\alpha \rangle):\alpha \in
{\underset\tilde {}\to S_{y,\zeta}} \right>$.  Now 
${\underset\tilde {}\to g_{y,\zeta}}$ is well defined, 
and actually can be computed if we use $A_\beta = \cup\{N^\alpha
_{\delta_\alpha}[{\underset\tilde {}\to G_{P_\zeta}}]:\alpha < \beta\}$.  So
by an assumption there is a suitable $P_\zeta$-name
${\underset\tilde {}\to C_\zeta}$ of a club of $\mu$.  But as $P_\zeta$
satisfies the $\mu$-c.c. \wilog \, ${\underset\tilde {}\to C_\zeta} =
C_\zeta$ so without loss of generality 
$C \subseteq \dbca_{\zeta \in w^*_y} C_\zeta$.  Now
choose $\alpha_1 < \alpha_2$ from $S_y \cap C$ and we choose by induction on
$\varepsilon \in w' = w^*_y \cup \{0,\ell g(\bar Q)\}$ a condition
$q_\varepsilon \in P_\varepsilon$ such that:

$$
\varepsilon_1 < \varepsilon \Rightarrow q_{\varepsilon_1} = q_\varepsilon
\restriction \varepsilon_1
$$

$$
q_\varepsilon \text{ is a bound to } \{p^{\alpha_1}_u \restriction
\varepsilon:i < \delta_{\alpha_1}\} \cup \{p^{\alpha_2}_i \restriction
\varepsilon:i < \delta_{\alpha_2}\}.
$$
\mn
For $\varepsilon = 0$ let $q_0 = \emptyset$, we have nothing to do really
if $\varepsilon$ is with no immediate predecessor in $w$, actually have
nothing to do $q_\varepsilon$ is $\cup\{q_{\varepsilon_1}:\varepsilon_1 <
\varepsilon,\varepsilon_2 \in w'\}$.  So let $\varepsilon = \varepsilon_1
+1,\varepsilon_1 \in w'$ and we use Definition \scite{3.1}.
\hfill$\square_{\scite{3.2}}$
\enddemo
\bn
\ub{\stag{3.2A} Comment on \scite{3.2}}:  1) There is ambiguity: 
$\hat{\Cal E}_0$ is as in the accessible case, \ub{but} this part works in
the other cases.  In particular, in Case A,B (in \cite{Sh:587}'s context)
if the length of $\bar a \in \hat{\Cal E}_0$ is $< \lambda$ (remember 
$\kappa = \lambda^+$), \ub{then} we
have $(< \lambda)$-completeness implies $\hat{\Cal E}_0$-completeness AND
in \scite{3.2} even $\bar a \in \hat{\Cal E}_0 \Rightarrow \ell g(\bar a) =
\omega$ is O.K.

In Case A on the $S_0 \subseteq S^\kappa_\lambda$ if $\ell g(\bar a) =
\lambda,a_\lambda \in S_0$ is O.K., too.  STILL can start with other variants
of completeness which is preserved.
\bigskip

\proclaim{\stag{3.3} Claim}  If $|Q| \le \kappa,{\Cal E}_0 
\subseteq \{\bar a:\bar a$
increasingly continuous $a_i \in [\mu^\kappa] < \kappa\}$ non-trivial possibly
just for one cofinality say $\aleph_0$, then $Q$ satisfies $\kappa^+$-pic.
\endproclaim
\bigskip

\demo{Proof} Trivial.  We get same sequence of condition.
\enddemo
\bn
\ub{\stag{3.4} Discussion}:  1) What is the use of pic?

In the forcing axioms instead ``$|Q| \le \kappa$" we can write "$Q$ satisfies
the $\kappa^+$-pic".  This strengthens the axioms.

In \cite{Sh:f} in some cases the length of the forcing is bounded
(there $\omega_2$) but here no need (as in \cite[Ch.VII,\S1]{Sh:f}).

This section applies to all cases in \cite{Sh:587} and its branches. \nl
2) Note that we can demand that the $P^\alpha_i$ satisfies some additional
requirements (in Definition \scite{3.1}) say $p^\alpha_{2i} = F_Q(\bar N
\restriction (2i+1),\bar p^\alpha \restriction (2i+1))$.
\newpage

\newpage

\head {\S4 Existence of non-free Whitehead (and Ext$(G,\Bbb Z) =
\{0\}$) abelian groups in successor of singulars} \endhead  \resetall
\bigskip

\proclaim{\stag{4.1} Claim}  Assume
\mr
\item "{$(a)$}"  $\lambda$ is strong limit singular, $\sigma = \text{ cf}
(\lambda) < \lambda,\kappa = \lambda^+ = 2^\lambda$
\sn
\item "{$(b)$}"  $S \subseteq \{\delta < \kappa:\text{cf}(\delta) =
\sigma\}$ is stationary
\sn
\item "{$(c)$}"  $S$ does not reflect or at least
\sn
\item "{$(c)^-$}"  $\bar A = \langle A_\delta:\delta \in S \rangle,
\text{otp}(A_\delta) = \sigma,\sup(A_\delta) = \delta$, \nl
$\bar A$ is $\lambda$-free, i.e. for every $\alpha^* < \kappa$ we can find
$\langle \alpha_\delta:\delta \in \alpha^* \cap S \rangle,\alpha_\delta <
\delta$ such that $\langle A_\delta \backslash \alpha_\delta:\delta \in S
\cap \alpha^* \rangle$ is a sequence of pairwise disjoint sets
\sn
\item "{$(d)$}"   $\langle G_i:i \le \sigma \rangle$ is a sequence of
abelian groups such that:
$$
\delta < \sigma \text{ limit } \Rightarrow G_\delta = \dbcu_{i < \delta} G_i
$$

$$
i < j \le \sigma \Rightarrow G_j/G_i \text{ free } \and G_i \subseteq G_j
$$
\sn
\item "{$(e)$}"  $G_\sigma / \dbcu_{i < \sigma} G_i \text{ is not Whitehead}$

$$
|G| < \lambda
$$

$$
\text{and } G_0 = \{0\}.
$$
\ermn
\ub{Then} \nl
1) There is a strongly $\kappa$-free abelian group of cardinality $\kappa$
which is not Whitehead, in fact $\Gamma(G) \subseteq S$. \nl
2) There is a strongly $\kappa$-free abelian group $G^*$ of cardinality
$\kappa$, HOM$(G^*,\Bbb Z) = \{0\}$, in fact $\Gamma(G) \subseteq S$ (in
fact the same abelian group can serve). \nl
\endproclaim
\bigskip

\remark{Remark}  We can replace ``not Whitehead" by: some $f \in \text{ HOM}
(\dbcu_{i < \sigma} G_i,\Bbb Z)$ cannot be extended to $f' \in \text{ HOM}
(G_\sigma,\Bbb Z)$.  This is used in part (2) which actually implies part (1).
\endremark
\bn
We first note:
\proclaim{\stag{4.2} Claim}  Assume
\mr
\item "{$(a)$}"  $\lambda$ strong limit singular, $\sigma = \text{ cf}
(\lambda) < \lambda,\kappa = 2^\lambda = \lambda^+$
\sn
\item "{$(b)$}"  $S \subseteq \{\delta < \kappa:\text{cf}(\delta) = \sigma
\text{ and } \lambda^\omega \text{ divides } \delta \text{ for simplicity}\}$
is stationary
\sn
\item "{$(c)$}"  $A_\delta \subseteq \delta = \sup(A_\delta)$, otp$(A_\delta)
 = \sigma,A_\delta = \{\alpha_{\delta,\zeta}:\zeta < \sigma\}$ increasing
with $\zeta$
\sn
\item "{$(d)$}" $h_0:\kappa \rightarrow \kappa$ and $h_1:\kappa \rightarrow
\sigma$ be such that \nl
$(\forall \alpha < \kappa)(\forall \zeta < \sigma)
(\forall \gamma \in (\alpha,\kappa))(\exists^\lambda \beta \in 
[\gamma,\gamma + \lambda])
(h_0(\beta) = \alpha \and h_1(\beta) = \zeta)$, \nl
$(\forall \alpha < \kappa) h_0(\alpha) \le \alpha$
\sn
\item "{$(e)$}"   Let $\lambda = \langle \lambda_\zeta:\zeta < 
\sigma \rangle$ increasingly continuous with limit $\lambda,\lambda_0 = 0$.
\ermn
\ub{Then} we can choose $\langle (g_\delta,\langle \beta^\delta_\zeta:\zeta
< \lambda \rangle):\delta \in S \rangle$ such that
\mr
\item "{$\bigodot_1$}"  $\langle \beta^\delta_\zeta:\zeta < \lambda \rangle$
is strictly increasing with limit $\delta$
\sn
\item "{$\bigodot_2$}"  for every $g:\kappa \rightarrow \kappa,B \in
[\kappa]^{< \lambda},g^2_\zeta:\kappa \rightarrow \lambda_{\zeta +1}$ for
$\zeta < \sigma$ \ub{there are} stationarily many $\delta \in S$ such that: 
{\roster
\itemitem{ $(i)$ }  $g_\delta = g \restriction B$
\sn
\itemitem{ $(ii)$ }  if $\dbcu_{\xi < \zeta} \lambda_\varepsilon \le \xi <
\lambda_\zeta$ then $h_0(\beta^\delta_\xi) = h_0(\beta^\delta
_{\lambda_\zeta}) =
\alpha_{\delta,\zeta},h_1(\beta^\delta_\xi) = h_1(\beta^\delta_
{\lambda_\zeta}) = \zeta,g^2_\zeta(\beta^\delta_\xi) = g^2_\zeta
(\beta^\delta_{\lambda_\zeta})$.
\endroster}
\endroster
\endproclaim
\bigskip

\remark{Remark}  Note that subtraction is meaningful, $(*)$ is quite
strong.
\endremark
\bigskip

\demo{Proof}  By the proofs of \scite{1.1}, \scite{1.2} (can use guessing
clubs by $\alpha_{\delta,\zeta}$'s, can demand that $\beta^\delta_{2 \zeta},
\beta^\delta_{2 \zeta +1} \in [\alpha_{\delta,\zeta},\alpha_{\delta,\zeta} +
\lambda)$.

For making the section selfcontained, we give a proof; note that this 
will prove \scite{1.2}, too.  Let ${\Cal H}(\kappa) = \dbcu_{\alpha < \kappa}
M_\alpha,M_\alpha$ is $\prec$-increasingly continuous, $\|M_\lambda\| = 
\lambda,\langle M_\beta:\beta \le \alpha \rangle \in M_{\alpha +1},
{}^{\sigma >}(M_{\alpha +1}) \subseteq M_{\alpha +1}$.
Let $\lambda = \dsize \sum_{i < \sigma} \lambda_i,\lambda_i$ increasingly
continuous, $\lambda_{i+1} > 2^{\lambda_i},\lambda_0 = 0,\lambda_1 > 
2^\sigma$.  Let $F_i:[\kappa]^{\lambda_i} \rightarrow [\kappa]^\kappa$ for
$i < \sigma$ (this towards $(*)_2$ of \scite{1.2}).  For
$\alpha < \lambda^+$, let $\alpha = \dbcu_{i < \sigma} a_{\alpha,i}$ such
that $|a_{\alpha,i}| \le \lambda_i$ and $a_{\alpha,i} \in M_{\alpha +1}$.  
Without loss of generality $\delta \in
S \Rightarrow \delta$ divisible by $\lambda^\omega$ (ordinal exponentiation).
For $\delta \in S$ let $\langle \beta^\delta_i:i < \sigma \rangle$ be
increasingly continuous with limit $\delta,\beta^\delta_i$ divisible by
$\lambda$ and $> 0$.  For $\delta \in S$ let $\langle b^\delta_i:i < \sigma
\rangle$ be such that: $b^\delta_i \subseteq \beta^\delta_i,|b^\delta_i| \le
\lambda_i,b^\delta_i$ is increasingly continuous in $i$ and $\delta = 
\dbcu_{i < \sigma} b^\delta_i$ (e.g. $b^\delta_i = \dbcu_{j_1,j_2 < i}
a_{\beta^\delta_{j_1,j_2}} \cup \lambda_i)$.  We further demand
$\beta^\delta_i \ge \lambda \Rightarrow \lambda_i \subseteq b^\delta_i 
\cap \lambda$.  Let $\langle
f^*_\alpha:\alpha < \lambda^+ \rangle$ list the two-place functions with
domain an ordinal $< \lambda^+$ and range $\subseteq \lambda^+$.  
Let $H$ be the set of functions $h$, Dom$(h) \in [\kappa]^{< \lambda}$,
Rang$(h) \subseteq \kappa$, so $|H| = \kappa$.  Let
$S = \cup \{S_h:h \in H\}$, with each $S_h$ stationary and
$\langle S_h:h \in H \rangle$ pairwise disjoint.  We now fix $h \in H$ 
and will choose $\bar \gamma^\delta = \langle \gamma^\delta_i:
i < \lambda \rangle$ for $\delta \in S_h$ such that clauses
$\bigodot_1 + \bigodot_2$ for our fix $h$ (and $\delta \in S_h$) hold, this
clearly suffices. \nl
Now for $\delta \in S_h$ and $i < \sigma$ and $g \in {}^{i+1}\sigma$ 
we can choose
$\zeta^\delta_{i,g,\varepsilon}$ (for $\varepsilon < \lambda_{i+1})$ 
such that:
\mr
\item "{$(A)$}"  $\langle \zeta^\delta_{i,g,\varepsilon}:\varepsilon < 
\lambda_i \rangle$ is a strictly increasing sequence of ordinals
\sn
\item "{$(B)$}"  $\beta^\delta_i < \zeta^\delta_{i,g,\varepsilon} <
\beta^\delta_{i+1}$, (can even demand
$\zeta^\delta_{i,j,\varepsilon} < \beta^\delta_i + \lambda$)
\sn
\item "{$(C)$}"  if $(\beta^\delta_i,\beta^\delta_{i+1}) \cap F_i
\{\zeta^\delta_{j,g \restriction (j+1),\varepsilon}:\varepsilon < \lambda_j
\text{ and } j < i\}$ has cardinality $\lambda$ then $\zeta^\delta_{i,g,
\varepsilon}$ belongs to it
\sn
\item "{$(D)$}"   for every $\alpha_1,\alpha_2 \in b^\delta_j$, the sequence 
$\langle \text{Min}\{\lambda_{g(i)},f^*_{\alpha_1}(\alpha_2,
\zeta^\delta_{i,g,\varepsilon}):\varepsilon < \lambda_{i+1}\}\rangle$ 
is constant i.e.: either
\endroster

$$
\varepsilon < \lambda_{i+1} \Rightarrow 
(\alpha_2,\zeta^\delta_{i,g,\varepsilon}) \notin \text{ Dom}(f^*_{\alpha_1})
$$

$$
\text{\ub{or} } \dsize \bigwedge_{\varepsilon < \lambda_{i+1}}
f^*_{\alpha_1}(\alpha_2,\zeta^\delta_{i,g,\varepsilon}) =
f^*_{\alpha_1}(\alpha_2,\zeta^\delta_{i,j,0})
$$

$$
\text{\ub{or} } \dsize \bigwedge_{\varepsilon < \lambda_j} f^*_{\alpha_1}
(\alpha_2,\zeta^\delta_{i,g,\varepsilon}) \ge \lambda_j.
$$
\mn
We can add $\langle f^*_{\alpha_1}(\alpha_2,\zeta^\delta_{i,g,\varepsilon}):
\varepsilon < \lambda_i \rangle$ is constant or strictly increasing (or all
undefined).
\mr
\item "{$(E)$}"  for some $j < \sigma,\{\varepsilon < \lambda_{i+1}:
\gamma^\delta_{i,g,\varepsilon} \in a_{\alpha,j}\}$ where \nl
$\alpha = \sup\{\gamma^\delta_{i,g,\varepsilon}:\varepsilon 
< \lambda_{i+1}\}$, (remember $\sigma \ne \lambda_{i+1}$ are regular).
\ermn
For each function $g \in {}^\sigma \sigma$ we try 
$\bar \gamma^{g,\delta} = \langle \gamma^{\delta,g}_\varepsilon:\varepsilon
< \lambda \rangle$ be: if $\lambda_i \le \varepsilon < \lambda_{i+1}$ then
$\gamma^{\delta,g}_\alpha = \gamma^\delta_{i,g,\varepsilon}$. \nl
Now for some $g$ it works.  \hfill$\square_{\scite{4.2}}$
\enddemo
\bigskip

\demo{Proof of \scite{4.1}}  1) We apply \scite{4.2} to the
$\langle A_\delta:\delta \in S \rangle$ from \scite{4.1}, and any reasonable
$h_0,h_1$. \nl
Let $\{t^{i,j}_\gamma + G_i:\gamma < \theta^{i,j}\}$ be a free basis of
$G^j/G^i$ for $i<j\le \sigma$.  If $i=0,j=\sigma$ we may omit the $i,j$, i.e. 
$t_\zeta = t^{0,\sigma}_\zeta$.  Let $\theta + \aleph_0 = |G| <a \lambda$
actually $\theta^{\zeta,\zeta+1} < \lambda_\zeta$ is enough, \wilog \,
$\theta < \lambda_1$ in \scite{4.2}.  Let $\beta^\delta_{\zeta,i} = 
\beta^\delta_{\xi(\zeta,i)}$ where $\xi(\zeta,i) = \dbcu_{\varepsilon < \zeta}
\lambda_\varepsilon + 1 + i$ for $\delta \in S,\zeta < \sigma,
i < \theta$.  \nl
Let $\beta_\delta(*) = \text{ Min}\{\beta:\beta \in \text{ Dom}(g_\delta),
g_\delta (\beta) \ne 0\}$, if well defined where $g_\delta$ is from 
\scite{4.2}.
\sn
Without loss of generality $\beta_\delta(*) \notin \{\beta^\delta_{\zeta,i}:
\zeta < \sigma,i < \theta\}$ (or omit $\lambda_\zeta,\beta^\delta_{\zeta,i}$ 
for $\zeta$ too small).  We define an abelian group $G^*$: it is 
generated by $\{x_\alpha:\alpha < \kappa\} \cup \{y^\delta_\gamma:\gamma < 
\theta^{0,\sigma} \text{ and } \delta \in S\}$ freely except the relations:
\mr
\item "{$(*)$}"   $\dsize \sum_{\gamma < \theta^{0,\sigma}} a_\gamma 
y^\delta_\gamma = \dsize \sum \bigl\{ b_{\zeta,\gamma}
(x_{\beta^\delta_{\zeta,\gamma}} - x_{\beta^\delta_{\lambda_\zeta}}):
\zeta < \sigma \text{ and } \gamma < \theta^{\zeta,\zeta +1} \bigr\}$ \nl
when $G \models \dsize \sum_{j < \theta} a_\gamma t_\gamma = 
\dsize \sum \bigl\{ b_{\zeta,\gamma} t^{\zeta,\zeta +1}_\gamma,\zeta <
\sigma \text{ and } \gamma < \theta^{\zeta,\zeta +1} \bigr\}$ where \nl
$a_\gamma,b_{\zeta,\gamma} \in \Bbb Z$ but all except finitely many are zero.
\ermn
There is a (unique) homomorphism $\bold g_\delta$ of $G$ into $H$ induced 
by $\bold g_\delta(t_\gamma) = g^\delta_\gamma$.  As usual it is an
embedding.  Let Rang $\bold g_\delta = G^{<\delta>}$. \nl
For $\alpha < \kappa$ let $G^*_\beta$ be the subgroup of $G^*$ generated by
$\{x_\alpha:\alpha < \delta\} \cup \{y^\delta_\gamma:\gamma <
\theta^{0,\sigma} \text{ and } \delta \in \beta \cap S\}$.  It can be
described like $G^*$.
\enddemo
\bn
\ub{Fact A}:  $G^*$ is strongly $\lambda$-free.
\bigskip

\demo{Proof}  For $\alpha^* < \beta^* < \kappa$, we can find $\langle
\alpha_\delta:\delta \in S \cap (\alpha^*,\beta^*] \rangle$ such that
$\langle A_\delta \backslash \alpha_\delta:\delta \in S \cap (\alpha^*,
\beta^*) \rangle$ are pairwise disjoint hence the sequence
$\langle \{\beta^\delta_{\zeta,i}:i < \theta,\zeta \in (\text{Min}
\{\xi < \sigma:\beta^\delta_{\zeta,0} > \alpha_\delta\},\sigma)\}:\delta \in
S \cap (\alpha^*,\beta^*] \rangle$ in a sequence of pairwise disjoint sets.
\nl
Let $\zeta_\delta = \text{ Min}\{\zeta:\beta^\delta_{\zeta,0} >
\alpha_\delta\} < \sigma$.  Now $G^*_{\beta^* +1}$ is generated as an
extension of $G^*_{\alpha^*+1}$ by $\{\bold g_\delta(t^{\zeta_\delta,\sigma}
_\gamma):\gamma < \theta^{\zeta_\delta,\sigma} \text{ and } \delta \in S \cap
(\alpha^*,\beta^*)\} \cup \{x_\alpha:\alpha \in (\alpha^*,\beta^*]$ and for 
no $\delta \in S \cap (\alpha^*,\beta^*]$ do we have $\alpha \in 
\{\beta^\delta_{\zeta,i}:i < \theta^{\zeta,\sigma} \text{ and } \zeta <
\zeta_\delta\}\}$; moreover $G^*_{\beta^*+1}$
is freely generated (as an extension of
$G^*_{\alpha^* +1}$).  So $G^*_{\beta^* +1}/G^*_{\alpha^*+1}$ is free, as
also $G^*_{\theta +1}$ is free we have shown Fact A.
\enddemo
\bn
\ub{Fact B}:  $G^*$ is not Whitehead.
\bigskip

\demo{Proof}  We choose by induction on $\alpha \le \kappa$, an abelian 
group $H_\alpha$ and a homomorphism $\bold h_\alpha:H_\alpha \rightarrow
G^*_\alpha = \langle \{x_\beta:\beta < \alpha\} \cup \{ y^\delta_\gamma:
\gamma < \theta,\delta \in S \cap \alpha\} \rangle_{G^*}$ 
increasingly continuous
in $\alpha$, with kernel $\Bbb Z,\bold h_0 =$ zero and $\bold k_\alpha:
G^*_\alpha \rightarrow H_\alpha$ not necessarily linear such that
$\bold h_\alpha \circ \bold k_\alpha = \text{ id}_{G^*_\alpha}$.  We identify
the set of members of $H_\alpha,G_\alpha,\Bbb Z$ with subsets of
$\lambda \times (1 + \alpha)$ such that $O_{H_\alpha} = O_{\Bbb Z} = 0$. \nl
Usually we have no freedom or no interesting freedom.  But we have for
$\alpha = \delta +1$, \nl
$\delta \in S$.  What we demand is ($G^{\langle \delta \rangle}$ - 
see before Fact A):
\mr
\item "{$(*)$}"  letting $H^{<\delta>} = \{x \in H_{\delta +1}:
\bold h_{\delta +1}(x) \in G^{<\delta>}\}$, if $s^* = g_\delta
(x_{\beta_\delta(*)}) \in \Bbb Z \backslash \{0\}$ ($g_\delta$ from 
\scite{4.2}(*)), \ub{then} there is no homomorphism $f_\delta:
G^{<\delta>} \rightarrow H^{<\delta>}$ such that
{\roster
\itemitem{ $(\alpha)$ }  $f_\delta(x_{\beta^\delta_{\zeta,i}}) - 
\bold k_\delta(x_{\beta^\delta_{\zeta,i}}) \in \Bbb Z$ is the same for all
$i \in (\dbcu_{\varepsilon < \zeta} \lambda_\varepsilon,\lambda_\zeta]$
\sn
\itemitem{ $(\beta)$ }  $\bold h_{\delta +1} \circ f_\delta = \text{ id}
_{G^{<\delta>}}$.
\endroster}
\ermn
[Why is this possible?  By non-Whiteheadness of $G^\sigma/\dbcu_{i < \sigma}
G^i$.] \nl
The rest should be clear.

\enddemo
\bigskip

\demo{Proof of \scite{4.1}(2)}  Of course, similar to that of
\scite{4.1}(1) but with some changes.
\mn
\ub{Step A}:  Without loss of generality there is a homomorphism $g^*$ from
$\dbcu_{i < \sigma} G^i$ to $\Bbb Z$ which cannot be extended to a
homormopshim from $G^\sigma$ to $\Bbb Z$. \nl
[Why?  Standard.]
\mn
\ub{Step B}:  During the construction of $G^*$, we choose $G^*_\alpha$ by
induction on $\alpha \le \kappa$, but if $g_\delta(0)$ from \scite{4.2} is 
a member of $G^*_\delta$ in $(*)$ we replace
$(x_{\beta^\delta_{\zeta,\gamma}} - x_{\beta^\delta_{\lambda_\zeta}})$ by
$\bigl( x_{\beta^\delta_{\lambda_\zeta + 1 + \gamma}} - 
x_{\beta^\delta_{\lambda_\zeta}} + g^*(t^{\zeta,\zeta +1}_\gamma) \times
g_\delta(0) \bigr)$. \nl
So if in the end $f:G^* \rightarrow \Bbb Z$ is a non-zero homomorphism, let
$x^* \in G^*$ be such that $f(x^*) \ne 0$ and $|f^*(x^*)|$ is minimal under
this, so \wilog \, it is 1.  So for some $\delta$

$$
f(g_\delta(0)) = 1_{\Bbb Z}
$$

$$
f(x_{\beta^\delta_{\lambda_\zeta + 1 + 1 + \gamma}}) = 
f(x_{\beta^\delta_{\lambda_\zeta}}) \text{ for } \gamma < \lambda_{\zeta +1}
- \lambda_\zeta
$$
\mn
(in fact for stationarily many).

So we get an easy contradiction. \hfill$\square_{\scite{4.1}}$
\enddemo
\bn
\ub{\stag{4.3} Comment}:  Assume
\mr
\item "{$\bigotimes_0$}"   $\lambda$ is strong limit singular, $\kappa =
2^\lambda = \lambda^+$ and cf$(\lambda) = \sigma$. 
\ermn
Assume
\mr
\item "{$\bigotimes_1$}"   $G$ is strongly $\kappa$-free not free of 
cardinality $\kappa$.
\ermn
\ub{Then} we can apply the analysis of \cite{Sh:161} (or see \cite{Sh:521}).  
So we have a $\kappa$-system ${\Cal S}$ and $\langle a^\ell_\eta:\eta \in
{\Cal S}_f,\ell < n({\Cal S}) \rangle,a^\ell_\eta$ countable sets as there.

Let $\langle G_i:i < \kappa \rangle$ be a filtration of $G$ so $G_i$ increases
continuously, $G_i$ a pure subgroup of $G,|G_i| < \kappa,G = 
\dbcu_{i < \kappa} G_i$, and for each $i,(\exists j > i)
(G_j/G_i \text{ non-free}) \rightarrow G_{i+1}/G_i$ non-free and
$(\exists j > i)(G_j/G_i$ not Whitehead $\rightarrow G_{i+1}/G_i$ not
Whitehead). \nl
Now
\mr
\item "{$(*)_1$}"  if $G$ is Whitehead, then $S_0 = \{\delta < \lambda^+:
\text{cf}(\delta) \ne \sigma,G_{\delta +1}/G_\delta$ not Whitehead$\}$ is not
stationary (otherwise by $\diamondsuit_{S_0}$ we are done).
\ermn
Assume further
\mr
\item "{$\bigotimes_2$}"  $G$ is Whitehead.
\ermn
Together
\mr
\item "{$(*)_2$}"  \wilog \, $\Gamma'(G) = \{\delta < \lambda^+:
G_{\delta +1}/G_\delta \text{ not Whitehead}\} \subseteq 
\{\delta < \kappa:\text{cf}(\delta) = \sigma\}$.
\ermn
Hence
\mr
\item "{$(*)_3$}"  \wilog \, $<>\ne \eta \in {\Cal S} \Rightarrow \eta(0)
\in \Gamma'(G)$ so cf$(\eta(0)) = \kappa$ and $\langle \delta \rangle \in S
\Rightarrow S^{[\delta]}$ analyze $G_{\delta +1}/G_\delta$
\ermn
hence (\cite{Sh:161}) \wilog \, for some $m(*) < n({\Cal S})$
\mr
\item "{$(*)_4$}"  $\eta \in {\Cal S}_f \Rightarrow \lambda^{\Cal S}
_{\eta \restriction m(*)} = \sigma$.
\ermn
Also \wilog
\mr
\item "{$(*)_5$}"  for all $\langle \delta \rangle \in {\Cal S}_f$ the
systems ${\Cal S}^{[<\delta >]}_f$ are isomorphic.
\ermn
This reduces us to essentially having clauses (d) + (e) of \scite{4.1}, more 
exactly from $S^{[<\delta>]}$ we can define $\langle G_i:i \le \sigma \rangle$
as in \scite{4.1}(d) except possibly (e) i.e. the Whiteheadness which is
guaranteed by
\mr
\item "{$\bigotimes^4_\kappa$}"  there is no nonfree Whitehead abelian group
of cardinality $< \kappa$.
\ermn
Generally also repeating (\cite{EkSh:505}) we can get that: \nl
\ub{Theorem} ($\lambda$ is strong limit singular $\kappa = 2^\lambda =
\lambda^+$, cf$(\lambda) = \sigma$ and $\bigotimes^4_\kappa$). \nl
The following are equivalent:
\mr
\widestnumber\item{$(b)(ii)$}
\item "{$(a)$}"  there is a strongly $\kappa$-free not free abelian group
of cardinality $\kappa$ which is Whitehead
\sn
\item "{$(b)(i)$}"  there is a strongly $\kappa$-free not free abelian group
of cardinality $\kappa$
\sn
\item "{${}(ii)$}"  there is $\langle G^i:i \le \sigma +1 \rangle$
increasingly continuous sequence of abelian groups such that
$i \le j \le \sigma +1 \and (i,j) \ne (\sigma,\sigma +1) 
\Rightarrow G^j/G^i$ free,$G^{\sigma +1}/G^\sigma$ not free, 
$|G^{\sigma +1}| \le \lambda$ and $G^{\sigma +1}/G^\sigma$ is Whitehead.
\ermn
[Why?  \ub{$(b) \Rightarrow (a)$} \,\, Analyze the group in $(i)$ getting
${\Cal S},\kappa$-system for freeness (instead of the earlier (in
\scite{4.1}) use $\bar A = \langle A_\delta:\delta \in S \rangle$ we now
represent $G$ and represent the $G^{\sigma +1}$ from $(ii)$ 
and puts them together, also with the
$\langle <\beta^\delta_t:i < \lambda>:\delta \in S \rangle$.  Better analyze
$G^{\sigma +1}/G^\sigma$ - like \cite{EkSh:505}.]
\mn
\ub{$(a) \Rightarrow (b)$}:  Essentially said earlier.] \nl
\newpage

\shlhetal

\newpage
    
REFERENCES.  
\bibliographystyle{lit-plain}
\bibliography{lista,listb,listx,listf,liste}

\enddocument

\bye